\documentclass[11pt]{article}
\usepackage{amscd}
\usepackage{amsfonts}
\usepackage{amsmath}
\usepackage{amssymb}
\usepackage{amsthm}
\usepackage{bbm}
\usepackage{CJK}
\usepackage{fancyhdr}
\usepackage{graphicx}
\usepackage{hyperref}
\usepackage{indentfirst}
\usepackage{latexsym}
\usepackage{mathrsfs}
\usepackage{xypic}
\usepackage{amsmath,amssymb,amsfonts,amsthm,fancyhdr}
\usepackage{epsfig,graphicx,picins,picinpar,subfigure}
\usepackage{pstricks}
\usepackage{amssymb}
\usepackage{xypic}
\usepackage[compress]{cite}

\newtheorem{theorem}{Theorem}[section]
\newtheorem{prop}[theorem]{Proposition}
\newtheorem{lemma}[theorem]{Lemma}
\newtheorem{coro}[theorem]{Corollary}
\newtheorem{prop-def}{Proposition-Definition}[section]

\theoremstyle{definition}
\newtheorem{defn}[theorem]{Definition}

\newtheorem{remark}[theorem]{Remark}
\newtheorem{exam}[theorem]{Example}

\usepackage[top=1in,bottom=1in,left=1.25in,right=1.25in]{geometry}
\def\<{\langle}
\def\>{\rangle}

\date{\today}
\begin{document}
\renewcommand{\baselinestretch}{1.2}
\renewcommand{\arraystretch}{1.0}
\title{\bf  Weighted $\mathcal{O}$-operators     on Hom-Lie triple systems}
\author{{\bf Wen Teng$^{1}$,  Jiulin Jin$^{2}$}\\
{\small 1. School of Mathematics and Statistics, Guizhou University of Finance and Economics} \\
{\small  Guiyang  550025, P. R. of China}\\
  {\small   Email: tengwen@mail.gufe.edu.cn  }\\
{\small 2.   School of Science,  Guiyang University,   Guiyang  550005,  China}\\
  {\small   Email:  j.l.jin@hotmail.com }}

 \maketitle
\begin{center}
\begin{minipage}{13.cm}

{\bf Abstract}
In this paper, we first introduce the notion of a  weighted $\mathcal{O}$-operator  on Hom-Lie triple systems  with respect to an action on another Hom-Lie triple system.  Next, we  construct a cohomology of weighted $\mathcal{O}$-operator   on Hom-Lie triple systems, we use the first cohomology group to classify linear deformations and we   investigate the obstruction class of an extendable
$n$-order  deformation. We end this paper by introducing a new algebraic structure, in connection with weighted $\mathcal{O}$-operators, called Hom-post-Lie triple system. We further
show that Hom-post-Lie triple systems can be derived from Hom-post-Lie algebras.

 \smallskip

 \smallskip

 \end{minipage}
 \end{center}
 \normalsize\vskip0.5cm

\section{Introduction}
\def\theequation{\arabic{section}. \arabic{equation}}
\setcounter{equation} {0}

Jacobson \cite{Jacobson1,Jacobson2}  first introduced the concept of Lie triple system. The current terminology was proposed by Yamguti \cite{Yamaguti}. In fact,   Lie triple system first appeared in Cartan's work \cite{Cartan}  on Riemannian Geometry, and noticed that the tangent space of symmetric space is Lie triple system.  Lie triple system has important applications in physics, such as elementary particle theory, the theory of
quantum mechanics and numerical analysis of differential equations.  More research on Lie triple systems  have
been developed, see \cite{Chen,Chtioui,Guo4,Lister, Hodge, Harris,Kubo,Lin,Sun,Teng,Teng1,Wu2,Wu22,Zhou,Zhang} and references cited therein.

 Hartwig et al. \cite{Hartwig} introduced the concept of Hom-Lie algebra as part of the research on the deformation of Witt and Virasoro algebras. Since then, Hom-type algebras have been widely studied because of their close relationship with discrete and deformed vector fields and differential calculus, see \cite{Cai, Hu, Sheng12}.
In particular, Sheng \cite{Sheng14}   proposed the representations  and cohomologies of Hom-Lie algebras.   Later, Ammar et al \cite{Ammar} considered the central extensions and deformations of Hom-Lie algebras from
cohomological points of view.   Das and his collaborators  \cite{Das22, Das23} explored the Nijenhuis operator and embedding tensor   on Hom-Lie algebras.  Nonabelian embedding tensors on Hom-Lie algebras were studied in \cite{Teng2}.

The concept of Rota-Baxter operator on associative algebra was introduced by Baxter \cite{B60} in his study of probability fluctuation theory.  Then the concept of an $\mathcal{O}$-operator (also called a relative Rota-Baxter operator) on a Lie algebra was  independently introduced by Kupershmidt \cite{K99}  to better understand the classical Yang-Baxter equation. Later, Mishra and Naolekar   \cite{Mishraa18} studied $\mathcal{O}$-operators on Hom-Lie algebras.  Chtioui et al. studied $\mathcal{O}$-operators on Lie triple systems in \cite{Chtioui}. Related results have been extended to the Hom version, such as relative Rota-Baxter operators on Hom-Lie-Yamaguti algebras \cite{Teng3} and  relative Rota-Baxter operators on Hom-Lie triple systems \cite{Li}.

Recently, due to the outstanding work of  \cite{Bai,Guo0,Das2, Wang}, scholars began to pay attention to Rota-Baxter operators with arbitrary weights.
After that, relative Rota-Baxter operators of nonzero weights on 3-Lie algebras \cite{Hou} and  relative Rota-Baxter operators of nonzero weights on Lie
triple systems \cite{Wu22} appear successively.    Furthermore,   Hou,   Sheng and Zhou \cite{Hou}  introduced the concept of 3-post-Lie algebras, and proved that
 a relative Rota-Baxter operator of nonzero weights on 3-Lie algebras induces a 3-post-Lie algebra.

Inspired by these works, the first purpose of this paper is to study weighted $\mathcal{O}$-operators   on Hom-Lie triple systems and associated structures. First, we introduce an action of a Hom-Lie triple system on another Hom-Lie triple system, which is different from the case of Lie triple system action \cite{Wu22}.
We introduce the notion of a $\kappa$-weighted $\mathcal{O}$-operator   from a Hom-Lie triple system to another Hom-Lie triple system with respect to an action, here $\kappa\in\mathbb{K}$ is a fixed scalar.  We further establish the cohomology and deformation theories for $\kappa$-weighted $\mathcal{O}$-operators  on Hom-Lie triple systems.

The second purpose is to investigate the Hom-post-Lie triple systems  and analyze
the relation with the aforementioned $\kappa$-weighted $\mathcal{O}$-operator. We introduce a new algebraic structure, which is called a Hom-post-Lie triple system. A Hom-post-Lie triple system naturally gives rise to a new Hom-Lie triple system and an action on the original Hom-Lie triple system such that the identity map is a $\kappa$-weighted $\mathcal{O}$-operator.  We show that
a $\kappa$-weighted $\mathcal{O}$-operator induces a Hom-post-Lie triple system structure naturally. Moreover, we
show that Hom-post-Lie triple systems can be derived from Hom-post-Lie algebras.

The paper is organized as follows.   In Section  \ref{sec: O-operator},  we  recall representations   of Lie triple systems and we  propose  an   action of a  Hom-Lie triple system on another Hom-Lie triple system. Then we  introduce the  concept of   a $\kappa$-weighted  $\mathcal{O}$-operator   from a   Hom-Lie triple system   to another   Hom-Lie triple system   with
respect to an action.
In Section \ref{sec: Cohomologies},   we  define the cohomologies of $\kappa$-weighted $\mathcal{O}$-operators   on  Hom-Lie triple systems.
 In Section \ref{sec: deformations},  we   study  linear deformations and higher order deformations of $\kappa$-weighted $\mathcal{O}$-operators   on Hom-Lie triple systems  via the cohomology theory established in the former
section.
 In Section \ref{sec: Hom-post-Lts},  we introduce the notion of  Hom-post-Lie triple system, which is the underlying
algebraic structure of $\kappa$-weighted $\mathcal{O}$-operators.
In Section  \ref{sec: Hom-post-Lie},  We show that Hom-Lie algebras, Hom-post-Lie algebras, Hom-Lie triple systems
and Hom-post-Lie triple systems are closely related.

Throughout this paper, $\mathbb{K}$ denotes a field of characteristic zero. All the    vector spaces  and
   (multi)linear maps are taken over $\mathbb{K}$.

\section{  Weighted  $\mathcal{O}$-operators   on Hom-Lie triple systems } \label{sec: O-operator}
\def\theequation{\arabic{section}.\arabic{equation}}
\setcounter{equation} {0}

In this section,  we  recall  concepts of    Hom-Lie triple systems from \cite{Yau} and \cite{Ma}.
Then, for   $\kappa\in \mathbb{K}$,  we introduce the  concept of a $\kappa$-weighted $\mathcal{O}$-operator  on  Hom-Lie triple system.

\begin{defn}  \cite{Yau}
(i) A Hom-Lie triple system (Hom-Lts) is a    triplet $(\mathfrak{g}, [-, -, -]_\mathfrak{g},\alpha_\mathfrak{g})$ in which $\mathfrak{g}$ is a vector space together with  a trilinear operation $[-, -, -]_\mathfrak{g}$
on $\mathfrak{g}$  and a linear map $\alpha_\mathfrak{g}: \mathfrak{g} \rightarrow  \mathfrak{g}$,  called the twisted map, satisfying $\alpha_\mathfrak{g}([x,y,z]_\mathfrak{g})=[\alpha_\mathfrak{g}(x),\alpha_\mathfrak{g}(y),\alpha_\mathfrak{g}(z)]_\mathfrak{g}$ such that
\begin{align}
&[x,y,z]_\mathfrak{g}+[y,x,z]_\mathfrak{g}=0,\label{2.1}\\
&[x,y,z]_\mathfrak{g}+[z,x,y]_\mathfrak{g}+[y,z,x]_\mathfrak{g}=0,\label{2.2}\\
 &[\alpha_\mathfrak{g}(a), \alpha_\mathfrak{g}(b), [x, y, z]_\mathfrak{g}]_\mathfrak{g}=[[a, b, x]_\mathfrak{g}, \alpha_\mathfrak{g}(y), \alpha_\mathfrak{g}(z)]_\mathfrak{g}+ [\alpha_\mathfrak{g}(x),  [a, b, y]_\mathfrak{g}, \alpha_\mathfrak{g}(z)]_\mathfrak{g}\nonumber\\
 &~~~~~~~~~~~~~~~~~~~~~~~~~~~~~~~~~~~~~~~~~~~~~~~~~~~~~~~~~~~~~~~+ [\alpha_\mathfrak{g}(x),  \alpha_\mathfrak{g}(y), [a, b, z]_\mathfrak{g}]_\mathfrak{g},\label{2.3}
\end{align}
where $ x, y, z, a, b\in \mathfrak{g}$.
In particular, the Hom-Lie triple system  $(\mathfrak{g}, [-, -, -]_\mathfrak{g},\alpha_\mathfrak{g})$  is said to be regular, if $\alpha_\mathfrak{g}$ is nondegenerate.

(ii) A homomorphism between two  Hom-Lie triple systems  $(\mathfrak{g}_1, [-, -, -]_1,\alpha_1)$ and $(\mathfrak{g}_2, [-, -,$ $-]_2,\alpha_2)$ is a linear map $\varphi: \mathfrak{g}_1\rightarrow \mathfrak{g}_2$ satisfying
$$\varphi(\alpha_1(x))=\alpha_2(\varphi(x)), ~~\varphi([x, y, z]_1)=[\varphi(x), \varphi(y),\varphi(z)]_2,~~\forall ~x,y,z\in \mathfrak{g}_1.$$
\end{defn}

Let $(\mathfrak{g}, [-,-]_\mathfrak{g},  \alpha_\mathfrak{g})$ be a     Hom-Lie   algebra, then
$( \mathfrak{g}, [-,-,-]_\mathfrak{g},  \alpha_\mathfrak{g})$ is a  Hom-Lie triple system,
where  $ [x,y,z]_{\mathfrak{g}}=[[x,y]_\mathfrak{g},\alpha_\mathfrak{g}(z)]_\mathfrak{g}, \forall x,y,z\in \mathfrak{g}.$

\begin{defn} \cite{Ma}
A representation of  a  Hom-Lie triple system  $(\mathfrak{g}, [-, -,-]_\mathfrak{g},\alpha_\mathfrak{g})$ on  a Hom-vector space $(V,\beta)$ is a bilinear map $\theta: \mathfrak{g}\times \mathfrak{g}\rightarrow \mathrm{End}(V)$, such that for all  $x,y,a,b\in \mathfrak{g}$
\begin{align}
& \theta(\alpha_\mathfrak{g}(x),\alpha_\mathfrak{g}(y))\circ\beta=\beta\circ\theta(x,y),\label{2.4}\\
& \theta(\alpha_\mathfrak{g}(a),\alpha_\mathfrak{g}(b))\theta(x,y)-\theta(\alpha_\mathfrak{g}(y),\alpha_\mathfrak{g}(b))\theta(x,a)-\theta(\alpha_\mathfrak{g}(x), [y,a,b]_\mathfrak{g})\circ\beta\nonumber\\
&+D(\alpha_\mathfrak{g}(y),\alpha_\mathfrak{g}(a))\theta(x,b)=0,\label{2.5}\\
& \theta(\alpha_\mathfrak{g}(a),\alpha_\mathfrak{g}(b))D(x,y)-D(\alpha_\mathfrak{g}(x),\alpha_\mathfrak{g}(y))\theta(a,b)+\theta([x,y,a]_\mathfrak{g},\alpha_\mathfrak{g}(b))\circ\beta\nonumber\\
&+\theta(\alpha_\mathfrak{g}(a),[x,y,b]_\mathfrak{g})\circ\beta=0,\label{2.6}\
\end{align}
 where $D(x,y)=\theta(y,x)-\theta(x,y)$. We also denote a representation of $\mathfrak{g}$ on $(V,\beta)$ by $(V,\beta; \theta)$.
\end{defn}

 Next, we introduce the concept of $\kappa$-weighted $\mathcal{O}$-operators   on Hom-Lie triple systems. We first propose the   action of a  Hom-Lie triple system on another Hom-Lie triple system.

\begin{defn}
Let $(\mathfrak{g}, [-, -, -]_\mathfrak{g},\alpha_\mathfrak{g})$ and $(\mathfrak{h}, [-, -, -]_\mathfrak{h},\alpha_\mathfrak{h})$  be two Hom-Lie triple systems.
Let $( \mathfrak{h},\alpha_\mathfrak{h}; \theta)$ be a representation of $\mathfrak{g}$. If $\theta$ meets the following equations
\begin{align}
& \theta(\alpha_\mathfrak{g}(x),\alpha_\mathfrak{g}(y))[u,v,w]_\mathfrak{h}\nonumber\\
&=[\theta(x,y)u,\alpha_\mathfrak{h}(v),\alpha_\mathfrak{h}(w)]_\mathfrak{h}+[\alpha_\mathfrak{h}(u),\theta(x,y)v,\alpha_\mathfrak{h}(w)]_\mathfrak{h}+[\alpha_\mathfrak{h}(u),\alpha_\mathfrak{h}(v),\theta(x,y)w]_\mathfrak{h},\label{2.7}\\
& \theta(\alpha_\mathfrak{g}(x),\alpha_\mathfrak{g}(y))[u,v,w]_\mathfrak{h}=[\alpha_\mathfrak{h}(u),\alpha_\mathfrak{h}(v),\theta(x,y)w]_\mathfrak{h}=0,\label{2.8}
\end{align}
 where $x,y\in \mathfrak{g},u,v,w\in \mathfrak{h}$,  then $\theta$ is called an action of $\mathfrak{g}$ on $\mathfrak{h}$.  We also denote an action of $\mathfrak{g}$ on $\mathfrak{h}$ by $(\mathfrak{h}, \alpha_\mathfrak{h}; \theta^\dag)$.
\end{defn}

By Eqs. \eqref{2.7} and \eqref{2.8},   we have the following equations:
\begin{align}
& D(\alpha_\mathfrak{g}(x),\alpha_\mathfrak{g}(y))[u,v,w]_\mathfrak{h}\nonumber\\
&=[D(x,y)u,\alpha_\mathfrak{h}(v),\alpha_\mathfrak{h}(w)]_\mathfrak{h}+[\alpha_\mathfrak{h}(u),D(x,y)v,\alpha_\mathfrak{h}(w)]_\mathfrak{h}+[\alpha_\mathfrak{h}(u),\alpha_\mathfrak{h}(v),D(x,y)w]_\mathfrak{h},\label{2.9}\\
& D(\alpha_\mathfrak{g}(x),\alpha_\mathfrak{g}(y))[u,v,w]_\mathfrak{h}=[\alpha_\mathfrak{h}(u),\alpha_\mathfrak{h}(v),D(x,y)w]_\mathfrak{h}=0.\label{2.10}
\end{align}

\begin{exam}
Let $(\mathfrak{g}, [-, -, -]_\mathfrak{g},\alpha_\mathfrak{g})$ be a Hom-Lie triple system.  We define linear map
$$\mathcal{R}: \mathfrak{g}\times \mathfrak{g} \rightarrow \mathrm{End}(\mathfrak{g}), (a,b)\mapsto (x\mapsto [x,a,b]_\mathfrak{g}),$$
with $\mathcal{L}(a,b)(x)=\mathcal{R}(b,a)x-\mathcal{R}(a,b)x=[a,b,x]_\mathfrak{g}$. Then,  $(\mathfrak{g},\alpha_\mathfrak{g}; \mathcal{R})$ is  a representation of the Hom-Lie triple system $\mathfrak{g}$, which is called the adjoint representation of  $\mathfrak{g}$.  Furthermore, if, for any $x_1,x_2,x_3\in \mathfrak{g}, $ $[[x_1,x_2,x_3]_\mathfrak{g},\alpha_\mathfrak{g}(a), \alpha_\mathfrak{g}(b)]_\mathfrak{g}=[\alpha_\mathfrak{g}(a), \alpha_\mathfrak{g}(b),[x_1,x_2,x_3]_\mathfrak{g}]_\mathfrak{g}=0$, then, by  Eq. \eqref{2.3},     $\mathcal{R}$ is an adjoint  action of $\mathfrak{g}$ on itself.
\end{exam}

Through direct calculation, an action can be described by semidirect product Hom-Lie triple system.

\begin{prop}
Let $\theta$ be an action of a Hom-Lie triple system $(\mathfrak{g}, [-, -, -]_\mathfrak{g},\alpha_\mathfrak{g})$ on another Hom-Lie triple system $(\mathfrak{h}, [-, -, -]_\mathfrak{h},\alpha_\mathfrak{h})$.
Then $ \mathfrak{g} \oplus  \mathfrak{h}$ is a Hom-Lie triple system under the following maps:
\begin{eqnarray*}
&&\alpha_\mathfrak{g}\oplus\alpha_\mathfrak{h}(x+u):=\alpha_\mathfrak{g}(x)+\alpha_\mathfrak{h}(u),\\
&&[x+u, y+v, z+w]_{\theta}:=[x, y,z]_\mathfrak{g}+D(x,y)u-\theta(x,z)v+\theta(y,z)u+\kappa[u,v,w]_\mathfrak{h},
\end{eqnarray*}
for all     $x,y,z\in  \mathfrak{g}$ and $u,v,w\in\mathfrak{h}$.  $(\mathfrak{g} \oplus  \mathfrak{h}, [-,-,-]_{\theta},\alpha_\mathfrak{g}\oplus\alpha_\mathfrak{h})$  is called
the  semidirect product  Hom-Lie triple system, and denoted by $ \mathfrak{g}\ltimes_{\theta} \mathfrak{h}$.
\end{prop}

\begin{defn}
(i) Let $\theta$ be an action of a Hom-Lie triple system $(\mathfrak{g}, [-, -, -]_\mathfrak{g},\alpha_\mathfrak{g})$ on another Hom-Lie triple system $(\mathfrak{h}, [-, -, -]_\mathfrak{h},\alpha_\mathfrak{h})$.
For $\kappa\in  \mathbb{K}$,  a linear map $\mathcal{A}:\mathfrak{h}\rightarrow \mathfrak{g}$ is called a  $\kappa$-weighted $\mathcal{O}$-operator    from a Hom-Lie triple system $\mathfrak{h}$ to another Hom-Lie triple system $\mathfrak{g}$ with
respect to an action $\theta$ if   the following equations hold:
\begin{align}
& \mathcal{A}\circ \alpha_\mathfrak{h}=\alpha_\mathfrak{g}\circ \mathcal{A},\label{2.11}\\
& [\mathcal{A}u,\mathcal{A}v,\mathcal{A}w]_\mathfrak{g}=\mathcal{A}(D(\mathcal{A}u,\mathcal{A}v)w-\theta(\mathcal{A}u,\mathcal{A}w)v+\theta(\mathcal{A}v,\mathcal{A}w)u+\kappa[u,v,w]_\mathfrak{h}),\label{2.12}
\end{align}
 where $u,v,w\in \mathfrak{h}$.

(ii) Let $\mathcal{A}_1$ and $\mathcal{A}_2$ be two $\kappa$-weighted $\mathcal{O}$-operators  from a Hom-Lie triple system $(\mathfrak{h}, [-, -, -]_\mathfrak{h},$ $\alpha_\mathfrak{h})$ to  another Hom-Lie triple system $(\mathfrak{g}, [-, -, -]_\mathfrak{g},\alpha_\mathfrak{g})$ with respect to an action $\theta$. A
homomorphism from $\mathcal{A}_1$ to $\mathcal{A}_2$
consists of two Hom-Lie triple system homomorphisms $\varphi_\mathfrak{h}: \mathfrak{h}\rightarrow\mathfrak{h}$ and $\varphi_\mathfrak{g}: \mathfrak{g}\rightarrow\mathfrak{g}$
such that
\begin{align}
& \varphi_\mathfrak{g}\circ\mathcal{A}_1=\mathcal{A}_2\circ\varphi_\mathfrak{h},\label{2.13}\\
& \varphi_\mathfrak{h}(\theta(x,y)u)=\theta(\varphi_\mathfrak{g}(x),\varphi_\mathfrak{g}(y))\varphi_\mathfrak{h}(u)), \forall x,y\in \mathfrak{g},u\in \mathfrak{h}. \label{2.14}
\end{align}
\end{defn}

By Eq. \eqref{2.14},   we have the following equation:
\begin{align}
& \varphi_\mathfrak{h}(D(x,y)u)=D(\varphi_\mathfrak{g}(x),\varphi_\mathfrak{g}(y))\varphi_\mathfrak{h}(u)), \forall x,y\in \mathfrak{g},u\in \mathfrak{h}. \label{2.15}
\end{align}

\begin{remark}
 (i)  Let  $\mathcal{A}:\mathfrak{h}\rightarrow \mathfrak{g}$  be  a $\kappa$-weighted  $\mathcal{O}$-operator   from a Hom-Lie triple system $(\mathfrak{h}, [-, -, -]_\mathfrak{h},\alpha_\mathfrak{h})$ to another Hom-Lie triple system $(\mathfrak{g}, [-, -, -]_\mathfrak{g},\alpha_\mathfrak{g})$ with
respect to an action $\theta$.  If the   bracket $[-,-, -]_\mathfrak{h}=0$, then   $\mathcal{A}$
becomes  an  $\mathcal{O}$-operator on Hom-Lie triple system.  See~\cite{Li} for more details.

(ii) When   $ \alpha_\mathfrak{g}= \mathrm{id}_\mathfrak{g}$,  $ \alpha_\mathfrak{h}= \mathrm{id}_\mathfrak{h}$, we get the notion of  $\kappa$-weighted  $\mathcal{O}$-operator  on Lie triple system.
 See~\cite{Wu22} for more details about $\kappa$-weighted $\mathcal{O}$-operator    on Lie triple system.
\end{remark}


\begin{exam} \label{exam:4-dimensional Hom-Lts}
Let $(\mathfrak{h}, [-, -, -]_\mathfrak{h},\alpha_\mathfrak{h})$  be a 4-dimensional Hom-Lie triple system with a basis  $\{\varepsilon_1,\varepsilon_2,\varepsilon_3,\varepsilon_4\}$  and the nonzero multiplication is given by
$$ [\varepsilon_1, \varepsilon_2, \varepsilon_1]_\mathfrak{h}=\varepsilon_4,~~ \alpha_\mathfrak{h}(\varepsilon_1)=\varepsilon_1,~~\alpha_\mathfrak{h}(\varepsilon_2)=-\varepsilon_2,~~\alpha_\mathfrak{h}(\varepsilon_3)=\varepsilon_3,~~\alpha_\mathfrak{h}(\varepsilon_4)=-\varepsilon_4.$$
 It is obvious that  $\mathcal{R}$ is an adjoint  action of $\mathfrak{h}$ on itself. Moreover,
$\mathcal{A}=\left(
                      \begin{array}{cccc}
                        0 & 0 & 0 & 0 \\
                        0 & 1 & 0 & 0\\
                        0 & 0 & 1 & 0\\
                        0 & 0 & 0 & 0\\
                      \end{array}
                    \right)$ is a $\kappa$-weighted $\mathcal{O}$-operator   from   $\mathfrak{h}$ to $\mathfrak{h}$ with
respect to an adjoint  action $\mathcal{R}$.
\end{exam}

\begin{prop}\label{prop:graph}
 A linear map $\mathcal{A}:\mathfrak{h}\rightarrow \mathfrak{g}$ is   a $\kappa$-weighted $\mathcal{O}$-operator   from a Hom-Lie triple system $(\mathfrak{h}, [-, -, -]_\mathfrak{h},\alpha_\mathfrak{h})$ to another Hom-Lie triple system $(\mathfrak{g}, [-, -, -]_\mathfrak{g},\alpha_\mathfrak{g})$ with
respect to an action $\theta$ if and only if the graph of $\mathcal{A}$
$$Gr(\mathcal{A})=\{\mathcal{A}u+u~|~u\in \mathfrak{h}\}$$
 is a subalgebra of the semidirect product  Hom-Lie triple system  $ \mathfrak{g}\ltimes_{\theta} \mathfrak{h}$.
\end{prop}

\begin{proof}
 Let  $\mathcal{A}:\mathfrak{h}\rightarrow \mathfrak{g}$ be a  linear map,  for any $u,v,w\in\mathfrak{h}$, we have
\begin{align*}
\alpha_\mathfrak{g}\oplus\alpha_\mathfrak{h}(\mathcal{A}u+u)=&\alpha_\mathfrak{g}\circ \mathcal{A}(u)+\alpha_\mathfrak{h}(u),\\
[\mathcal{A}u+u,\mathcal{A}v+v,\mathcal{A}w+w]_{\theta}=&[\mathcal{A}u,\mathcal{A}v,\mathcal{A}w]_\mathfrak{g}+D(\mathcal{A}u,\mathcal{A}v)w-\theta(\mathcal{A}u,\mathcal{A}w)v\\
&+\theta(\mathcal{A}v,\mathcal{A}w)u+\kappa[u,v,w]_\mathfrak{h},
\end{align*}
which implies that the graph  $Gr( \mathcal{A})$ is a subalgebra of the semidirect product  Hom-Lie triple system  $ \mathfrak{g}\ltimes_{\theta} \mathfrak{h}$
if and only if $ \mathcal{A}$ satisfies Eqs.  \eqref{2.11} and \eqref{2.12}, which means that $ \mathcal{A}$ is a $\kappa$-weighted $\mathcal{O}$-operator.
\end{proof}

Since $Gr(\mathcal{A})$ is isomorphic to $\mathfrak{h}$ as a vector space. Define a trilinear operation on $\mathfrak{h}$ by
\begin{align*}
&\{u, v,  w\}_\mathcal{A}=D(\mathcal{A}u,\mathcal{A}v)w-\theta(\mathcal{A}u,\mathcal{A}w)v+\theta(\mathcal{A}v,\mathcal{A}w)u+\kappa[u,v,w]_\mathfrak{h}, ~\forall~ u,v,w\in\mathfrak{h}.
\end{align*}
By Proposition \ref{prop:graph}, we get $(\mathfrak{h}, \{-, -, -\}_\mathcal{A}, \alpha_\mathfrak{h})$ is a Hom-Lie triple system.

In \cite{Hou20}, Hou, Ma and Chen   introduced the notion of Nijenhuis operator through the 2-order  deformation on Hom-Lie triple system $(\mathfrak{g}, [-, -, -]_\mathfrak{g},\alpha_\mathfrak{g})$.
 More precisely, a linear map $N:\mathfrak{g}\rightarrow\mathfrak{g}$ is
called a  Nijenhuis operator on $\mathfrak{g}$ if for all $x,y,z\in \mathfrak{g}$, the following equations hold:
\begin{align*}
 N\circ \alpha_\mathfrak{g}=&\alpha_\mathfrak{g}\circ N,\\
 [Nx,Ny,Nz]_\mathfrak{g}=&N([x,Ny,Nz]_\mathfrak{g}+[Nx,y,Nz]_\mathfrak{g}+[Nx,Ny,z]_\mathfrak{g})\\
 &-N^2([Nx,y,z]_\mathfrak{g}+[x,Ny,z]_\mathfrak{g}+[x,y,Nz]_\mathfrak{g})+N^3[x,y,z]_\mathfrak{g}.
\end{align*}

Through direct verification, we have the relationship between $\kappa$-weighted $\mathcal{O}$-operators   and Nijenhuis operators.

\begin{prop}
 A linear map $\mathcal{A}:\mathfrak{h}\rightarrow \mathfrak{g}$ is  a $\kappa$-weighted $\mathcal{O}$-operator   from a Hom-Lie triple system $(\mathfrak{h}, [-, -, -]_\mathfrak{h},\alpha_\mathfrak{h})$ to another Hom-Lie triple system $(\mathfrak{g}, [-, -, -]_\mathfrak{g},\alpha_\mathfrak{g})$ with
respect to an action $\theta$ if and only if
$$ N_\mathcal{A}=\left(
        \begin{array}{cc}
          \mathrm{id}&  \mathcal{A}\\
           0 & 0  \\
        \end{array}
      \right)$$
 is a Nijenhuis operator on the semidirect product  Hom-Lie triple system  $ \mathfrak{g}\ltimes_{\theta} \mathfrak{h}$.
\end{prop}

\section{  Cohomologies of    weighted $\mathcal{O}$-operators   on Hom-Lie triple systems} \label{sec: Cohomologies}
\def\theequation{\arabic{section}.\arabic{equation}}
\setcounter{equation} {0}

In this section, we construct a representation of the    Hom-Lie triple system $( \mathfrak{h}, \{-, -, -\}_\mathcal{A},$ $\alpha_\mathfrak{h})$  on the Hom-vector space $(\mathfrak{g},  \alpha_\mathfrak{g})$
from a $\kappa$-weighted $\mathcal{O}$-operator $\mathcal{A} $    and define the cohomologies of
$\kappa$-weighted $\mathcal{O}$-operator    on  the Hom-Lie triple system.

\begin{lemma}
  Let  $\mathcal{A}:\mathfrak{h}\rightarrow \mathfrak{g}$  be  a $\kappa$-weighted  $\mathcal{O}$-operator   from a Hom-Lie triple system $(\mathfrak{h}, [-, -, -]_\mathfrak{h},\alpha_\mathfrak{h})$ to another Hom-Lie triple system $(\mathfrak{g}, [-, -, -]_\mathfrak{g},\alpha_\mathfrak{g})$ with
respect to an action $\theta$. For any $u,v\in\mathfrak{h}$, $x\in \mathfrak{g}$, define $\theta_\mathcal{A}: \mathfrak{h} \otimes \mathfrak{h}\rightarrow \mathrm{End}(\mathfrak{g})$ by
\begin{align}
\theta_\mathcal{A}(u,v)(x)=[x,\mathcal{A}u,\mathcal{A}v]_\mathfrak{g}+\mathcal{A}(\theta(x,\mathcal{A}v)u-D(x,\mathcal{A}u)v), \label{4.1}
\end{align}
then, $(\mathfrak{g},  \alpha_\mathfrak{g}; \theta_\mathcal{A})$ is a representation of the   Hom-Lie triple system $( \mathfrak{h}, \{ -,-,-\}_\mathcal{A},\alpha_\mathfrak{h})$.
\end{lemma}

\begin{proof}
For any $u,v,s,t\in \mathfrak{h}, x\in \mathfrak{g}$,   we have
\begin{align}
&D_\mathcal{A}(u,v)(x)\nonumber\\
=&\theta_\mathcal{A}(v,u)(x)-\theta_\mathcal{A}(u,v)(x)\nonumber\\
=&[x,\mathcal{A}v,\mathcal{A}u]_\mathfrak{g}+\mathcal{A}(\theta(x,\mathcal{A}u)v-D(x,\mathcal{A}v)u)-[x,\mathcal{A}u,\mathcal{A}v]_\mathfrak{g}-\mathcal{A}(\theta(x,\mathcal{A}v)u-D(x,\mathcal{A}u)v)\nonumber\\
=&[\mathcal{A}u,\mathcal{A}v,x]_\mathfrak{g}+\mathcal{A}(\theta(\mathcal{A}u, x)v-\theta(\mathcal{A}v, x)u). \label{4.2}
\end{align}
Further, we obtain that
\begin{small}
\begin{align*}
&\theta_\mathcal{A}(\alpha_\mathfrak{h}(u),\alpha_\mathfrak{h}(v))\alpha_\mathfrak{g}(x)\\
=&[\alpha_\mathfrak{g}(x),\mathcal{A}\alpha_\mathfrak{h}(u),\mathcal{A}\alpha_\mathfrak{h}(v)]_\mathfrak{g}+\mathcal{A}(\theta(\alpha_\mathfrak{g}(x),\mathcal{A}\alpha_\mathfrak{h}(v))\alpha_\mathfrak{h}(u)-D(\alpha_\mathfrak{g}(x),\mathcal{A}\alpha_\mathfrak{h}(u))\alpha_\mathfrak{h}(v))\\
=&[\alpha_\mathfrak{g}(x),\alpha_\mathfrak{g}(\mathcal{A}u),\alpha_\mathfrak{g}(\mathcal{A}v)]_\mathfrak{g}+\mathcal{A}(\theta(\alpha_\mathfrak{g}(x),\alpha_\mathfrak{g}(\mathcal{A}v))\alpha_\mathfrak{h}(u)-D(\alpha_\mathfrak{g}(x),\alpha_\mathfrak{g}(\mathcal{A}u))\alpha_\mathfrak{h}(v))\\
=&\alpha_\mathfrak{g}([x,\mathcal{A}u,\mathcal{A}v]_\mathfrak{g}+\mathcal{A}(\theta(x,\mathcal{A}v)u-D(x, \mathcal{A}u))v)\\
=&\alpha_\mathfrak{g}(\theta_\mathcal{A}(u,v)x),
\end{align*}
\end{small}
\begin{small}
\begin{align*}
& \theta_\mathcal{A}(\alpha_\mathfrak{h}(u),\alpha_\mathfrak{h}(v))\theta_\mathcal{A}(s,t)x-\theta_\mathcal{A}(\alpha_\mathfrak{h}(t),\alpha_\mathfrak{h}(v))\theta(s,u)x-\theta_\mathcal{A}(\alpha_\mathfrak{h}(s), [t,u,v]_\mathcal{A})\alpha_\mathfrak{g}(x) \\
&+D_\mathcal{A}(\alpha_\mathfrak{h}(t),\alpha_\mathfrak{h}(u))\theta_\mathcal{A}(s,v)x\\
=&[[x,\mathcal{A}s,\mathcal{A}t]_\mathfrak{g},\alpha_\mathfrak{g}(\mathcal{A}u),\alpha_\mathfrak{g}(\mathcal{A}v)]_\mathfrak{g}+\mathcal{A}(\theta([x,\mathcal{A}s,\mathcal{A}t]_\mathfrak{g},\alpha_\mathfrak{g}(\mathcal{A}v))\alpha_\mathfrak{h}(u)-D([x,\mathcal{A}s,\mathcal{A}t]_\mathfrak{g},\alpha_\mathfrak{g}(\mathcal{A}u))\alpha_\mathfrak{h}(v))\\
&+[\mathcal{A}(\theta(x,\mathcal{A}t)s),\alpha_\mathfrak{g}(\mathcal{A}u),\alpha_\mathfrak{g}(\mathcal{A}v)]_\mathfrak{g}+\mathcal{A}(\theta(\mathcal{A}(\theta(x,\mathcal{A}t)s),\alpha_\mathfrak{g}(\mathcal{A}v))\alpha_\mathfrak{h}(u)-D(\mathcal{A}(\theta(x,\mathcal{A}t)s),\alpha_\mathfrak{g}(\mathcal{A}u))\alpha_\mathfrak{h}(v))\\
&-[\mathcal{A}(D(x,\mathcal{A}s)t),\alpha_\mathfrak{g}(\mathcal{A}u),\alpha_\mathfrak{g}(\mathcal{A}v)]_\mathfrak{g}-\mathcal{A}(\theta(\mathcal{A}(D(x,\mathcal{A}s)t),\alpha_\mathfrak{g}(\mathcal{A}v))\alpha_\mathfrak{h}(u)-D(\mathcal{A}(D(x,\mathcal{A}s)t),\alpha_\mathfrak{g}(\mathcal{A}u))\alpha_\mathfrak{h}(v))\\
&-[[x,\mathcal{A}s,\mathcal{A}u]_\mathfrak{g},\alpha_\mathfrak{g}(\mathcal{A}t),\alpha_\mathfrak{g}(\mathcal{A}v)]_\mathfrak{g}-\mathcal{A}(\theta([x,\mathcal{A}s,\mathcal{A}u]_\mathfrak{g},\alpha_\mathfrak{g}(\mathcal{A}v))\alpha_\mathfrak{h}(t)-D([x,\mathcal{A}s,\mathcal{A}u]_\mathfrak{g},\alpha_\mathfrak{g}(\mathcal{A}t))\alpha_\mathfrak{h}(v))\\
&-[\mathcal{A}(\theta(x,\mathcal{A}u)s),\alpha_\mathfrak{g}(\mathcal{A}t),\alpha_\mathfrak{g}(\mathcal{A}v)]_\mathfrak{g}-\mathcal{A}(\theta(\mathcal{A}(\theta(x,\mathcal{A}u)s),\alpha_\mathfrak{g}(\mathcal{A}v))\alpha_\mathfrak{h}(t)-D(\mathcal{A}(\theta(x,\mathcal{A}u)s),\alpha_\mathfrak{g}(\mathcal{A}t))\alpha_\mathfrak{h}(v))\\
&+[\mathcal{A}(D(x,\mathcal{A}s)u),\alpha_\mathfrak{g}(\mathcal{A}t),\alpha_\mathfrak{g}(\mathcal{A}v)]_\mathfrak{g}+\mathcal{A}(\theta(\mathcal{A}(D(x,\mathcal{A}s)u),\alpha_\mathfrak{g}(\mathcal{A}v))\alpha_\mathfrak{h}(t)-D(\mathcal{A}(D(x,\mathcal{A}s)u),\alpha_\mathfrak{g}(\mathcal{A}t))\alpha_\mathfrak{h}(v))\\
&-[\alpha_\mathfrak{g}(x),\alpha_\mathfrak{g}(\mathcal{A}s),\mathcal{A}[t,u,v]_\mathcal{A}]_\mathfrak{g}-\mathcal{A}(\theta(\alpha_\mathfrak{g}(x),\mathcal{A}[t,u,v]_\mathcal{A})\alpha_\mathfrak{h}(s)-D(\alpha_\mathfrak{g}(x),\alpha_\mathfrak{g}(\mathcal{A}s))[t,u,v]_\mathcal{A})\\
&+[\alpha_\mathfrak{g}(\mathcal{A}t),\alpha_\mathfrak{g}(\mathcal{A}u),[x,\mathcal{A}s,\mathcal{A}v]_\mathfrak{g}]_\mathfrak{g}+\mathcal{A}(\theta(\alpha_\mathfrak{g}(\mathcal{A}t), [x,\mathcal{A}s,\mathcal{A}v]_\mathfrak{g})\alpha_\mathfrak{h}(u)-\theta(\alpha_\mathfrak{g}(\mathcal{A}u), [x,\mathcal{A}s,\mathcal{A}v]_\mathfrak{g})\alpha_\mathfrak{h}(t))\\
&+[\alpha_\mathfrak{g}(\mathcal{A}t),\alpha_\mathfrak{g}(\mathcal{A}u),\mathcal{A}(\theta(x,\mathcal{A}v)s)]_\mathfrak{g}+\mathcal{A}(\theta(\alpha_\mathfrak{g}(\mathcal{A}t), \mathcal{A}(\theta(x,\mathcal{A}v)s))\alpha_\mathfrak{h}(u)-\theta(\alpha_\mathfrak{g}(\mathcal{A}u), \mathcal{A}(\theta(x,\mathcal{A}v)s))\alpha_\mathfrak{h}(t))\\
&-[\alpha_\mathfrak{g}(\mathcal{A}t),\alpha_\mathfrak{g}(\mathcal{A}u),\mathcal{A}(D(x,\mathcal{A}s)v)]_\mathfrak{g}-\mathcal{A}(\theta(\alpha_\mathfrak{g}(\mathcal{A}t), \mathcal{A}(D(x,\mathcal{A}s)v))\alpha_\mathfrak{h}(u)-\theta(\alpha_\mathfrak{g}(\mathcal{A}u), \mathcal{A}(D(x,\mathcal{A}s)v))\alpha_\mathfrak{h}(t))\\
=&0.
\end{align*}
\end{small}
Similarly, we also have
\begin{align*}
&\theta_\mathcal{A}(\alpha_\mathfrak{h}(u),\alpha_\mathfrak{h}(v))D_\mathcal{A}(s,t)x-D_\mathcal{A}(\alpha_\mathfrak{h}(s),\alpha_\mathfrak{h}(t))\theta_\mathcal{A}(u,v)x+\theta_\mathcal{A}([s,t,u]_\mathcal{A},\alpha_\mathfrak{h}(v))\alpha_\mathfrak{g}(x) \\
&+\theta_\mathcal{A}(\alpha_\mathfrak{g}(u),[s,t,v]_\mathcal{A})\alpha_\mathfrak{g}(x)=0.
\end{align*}
Hence,  $(\mathfrak{g},  \alpha_\mathfrak{g}; \theta_\mathcal{A})$ is a representation of   $( \mathfrak{h}, \{ -,-,-\}_\mathcal{A},\alpha_\mathfrak{h})$.
\end{proof}

Let $(\mathfrak{g},  \alpha_\mathfrak{g}; \theta_\mathcal{A})$  be a representation of a  Hom-Lie triple system $( \mathfrak{h}, \{ -,-,-\}_\mathcal{A},\alpha_\mathfrak{h})$.
Denote the $(2n+1)$-cochains   of $ \mathfrak{h}$ with coefficients in representation $( \mathfrak{g}, \alpha_\mathfrak{g};  \theta_\mathcal{A})$   by
\begin{align*}
&\mathcal{C}_{\mathrm{HLts}}^{2n+1}( \mathfrak{h},  \mathfrak{g}):=\big\{f\in \mathrm{Hom}( \mathfrak{h}^{\otimes2n+1}, \mathfrak{g}) ~|~ \alpha_\mathfrak{g}(f(v_1,\cdots, v_{2n+1})=f(\alpha_\mathfrak{h}(v_1),\cdots, \alpha_\mathfrak{h}(v_{2n+1})),\\
&~f(v_1,\cdots, v_{2n-2},u,v,w)+f(v_1,\cdots, v_{2n-2},v,u,w)=0,~\circlearrowright_{u,v,w}f(v_1,\cdots, v_{2n-2},u,v,w)=0  \big\}.
\end{align*}


For $n\geq 1,$ let $\delta: \mathcal{C}_{\mathrm{HLts}}^{2n-1}(\mathfrak{h},  \mathfrak{g})\rightarrow \mathcal{C}_{\mathrm{HLts}}^{2n+1}( \mathfrak{h},  \mathfrak{g})$ be the corresponding  coboundary operator of the descent Hom-Lie triple system $( \mathfrak{h}, \{ -,-,-\}_\mathcal{A},\alpha_\mathfrak{h})$ with coefficients in the representation $(\mathfrak{g},  \alpha_\mathfrak{g}; \theta_\mathcal{A})$,   More precisely,  for $v_1,\cdots,v_{2n+1}\in  \mathfrak{h}$ and $f\in \mathcal{C}_{\mathrm{HLts}}^{2n-1}( \mathfrak{h},  \mathfrak{g})$, as
\begin{align*}
&\delta f(v_1,\cdots, v_{2n+1})\\
=&\theta_\mathcal{A}(\alpha^{n-1}_{\mathfrak{h}}(v_{2n}),\alpha^{n-1}_{\mathfrak{h}}(v_{2n+1}))f(v_1,\cdots, v_{2n-1})-\theta_\mathcal{A}(\alpha^{n-1}_{\mathfrak{h}}(v_{2n-1}),\alpha^{n-1}_{\mathfrak{h}}(v_{2n+1}))f(v_1,\cdots, v_{2n-2},v_{2n})\\
&+\sum_{i=1}^n(-1)^{i+n}D_\mathcal{A}(\alpha^{n-1}_{\mathfrak{h}}(v_{2i-1}),\alpha^{n-1}_{\mathfrak{h}}(v_{2i}))f(v_1,\cdots, v_{2i-2},v_{2i+1},\cdots, v_{2n+1})\\
&+\sum_{i=1}^n\sum_{j=2i+1}^{2n+1}(-1)^{i+n+1}f(\alpha_{\mathfrak{h}}(v_1),\cdots, \alpha_{\mathfrak{h}}(v_{2i-2}),\alpha_{\mathfrak{h}}(v_{2i+1}),\cdots, \{v_{2i-1},v_{2i},v_{j}\}_\mathcal{A},\cdots,\alpha_{\mathfrak{h}}(v_{2n+1})).
\end{align*}
So $\delta\circ\delta=0.$
One can refer to \cite{Ma} for more information about Hom-Lie triple systems and cohomology theory.

In particular, for $f\in \mathcal{C}_{\mathrm{HLts}}^{1}( \mathfrak{h}, \mathfrak{g})$, $f$ is a 1-cocycle  on  $( \mathfrak{h}, \{ -,-,-\}_\mathcal{A},\alpha_\mathfrak{h})$ with coefficients in   $(\mathfrak{g},  \alpha_\mathfrak{g}; \theta_\mathcal{A})$ if $\delta f (v_1,v_2,v_3)=0,$ i.e.,
$$\theta_\mathcal{A}(v_2,v_3)f(v_1)-\theta_\mathcal{A}(v_1,v_3)f(v_2)+D_\mathcal{A}(v_1,v_2)f(v_3)-f(\{v_{1},v_{2},v_{3}\}_\mathcal{A})=0$$
and
$$\alpha_\mathfrak{g}(f(v_1))-f(\alpha_\mathfrak{h}(v_1))=0.$$

\begin{prop}\label{prop:1-cocycle}
 Let  $\mathcal{A}:\mathfrak{h}\rightarrow \mathfrak{g}$  be  a $\kappa$-weighted  $\mathcal{O}$-operator    from a regular Hom-Lie triple system $(\mathfrak{h}, [-, -, -]_\mathfrak{h},\alpha_\mathfrak{h})$ to another regular Hom-Lie triple system $(\mathfrak{g}, [-, -, -]_\mathfrak{g},\alpha_\mathfrak{g})$ with
respect to an action $\theta$.  If  there exist  two elements $a,b\in  \mathfrak{g}$ such that $\alpha_\mathfrak{g}(a)=a,\alpha_\mathfrak{g}(b)=b$, we define $\Im(a,b): \mathfrak{h}\rightarrow  \mathfrak{g}$ by
$$\Im(a,b)v=\mathcal{A}(D(a,b)\alpha_\mathfrak{h}^{-1}(v))-[a,b,\mathcal{A}\alpha_\mathfrak{h}^{-1}(v)]_\mathfrak{g},$$
for any $v\in  \mathfrak{h}.$ Then $\Im(a,b)$ is a 1-cocycle on  $( \mathfrak{h}, \{ -,-,-\}_\mathcal{A},\alpha_\mathfrak{h})$ with coefficients in   $(\mathfrak{g},  \alpha_\mathfrak{g}; \theta_\mathcal{A})$.
\end{prop}

\begin{proof}
For any $v_1,v_2,v_3\in \mathfrak{h}$, first, obviously $\alpha_\mathfrak{g}(\Im(a,b)(v_1))-\Im(a,b)(\alpha_\mathfrak{h}(v_1))=0$.   Next,  by Eqs. \eqref{3.6}, \eqref{4.1} and \eqref{4.2}, we have
\begin{align*}
&\delta \Im(a,b) (v_1,v_2,v_3) \\
=& \theta_\mathcal{A}(v_2,v_3)\Im(a,b)v_1-\theta_\mathcal{A}(v_1,v_3)\Im(a,b)v_2+D_\mathcal{A}(v_1,v_2)\Im(a,b)v_3-\Im(a,b)\{v_{1},v_{2},v_{3}\}_\mathcal{A}\\
=& [\mathcal{A}(D(a,b)\alpha_\mathfrak{h}^{-1}(v_1)),\mathcal{A}v_2,\mathcal{A}v_3]_\mathfrak{g}+\mathcal{A}(\theta(\mathcal{A}(D(a,b)\alpha_\mathfrak{h}^{-1}(v_1)),\mathcal{A}v_3)v_2-D(\mathcal{A}(D(a,b)\alpha_\mathfrak{h}^{-1}(v_1)),\mathcal{A}v_2)v_3)\\
&-[[a,b,\mathcal{A}\alpha_\mathfrak{h}^{-1}(v_1)]_\mathfrak{g},\mathcal{A}v_2,\mathcal{A}v_3]_\mathfrak{g}-\mathcal{A}(\theta([a,b,\mathcal{A}\alpha_\mathfrak{h}^{-1}(v_1)]_\mathfrak{g},\mathcal{A}v_3)v_2-D([a,b,\mathcal{A}\alpha_\mathfrak{h}^{-1}(v_1)]_\mathfrak{g},\mathcal{A}v_2)v_3)\\
&-[\mathcal{A}(D(a,b)\alpha_\mathfrak{h}^{-1}(v_2)),\mathcal{A}v_1,\mathcal{A}v_3]_\mathfrak{g}-\mathcal{A}(\theta(\mathcal{A}(D(a,b)\alpha_\mathfrak{h}^{-1}(v_2)),\mathcal{A}v_3)v_1-D(\mathcal{A}(D(a,b)\alpha_\mathfrak{h}^{-1}(v_2)),\mathcal{A}v_1)v_3)\\
&+[[a,b,\mathcal{A}\alpha_\mathfrak{h}^{-1}(v_2)]_\mathfrak{g},\mathcal{A}v_1,\mathcal{A}v_3]_\mathfrak{g}+\mathcal{A}(\theta([a,b,\mathcal{A}\alpha_\mathfrak{h}^{-1}(v_2)]_\mathfrak{g},\mathcal{A}v_3)v_1-D([a,b,\mathcal{A}\alpha_\mathfrak{h}^{-1}(v_2)]_\mathfrak{g},\mathcal{A}v_1)v_3)\\
&+[\mathcal{A}v_1,\mathcal{A}v_2,\mathcal{A}(D(a,b)\alpha_\mathfrak{h}^{-1}(v_3))]_\mathfrak{g}+\mathcal{A}(\theta(\mathcal{A}v_1, \mathcal{A}(D(a,b)\alpha_\mathfrak{h}^{-1}(v_3)))v_2-\theta(\mathcal{A}v_2, \mathcal{A}(D(a,b)\alpha_\mathfrak{h}^{-1}(v_3)))v_1)\\
&-[\mathcal{A}v_1,\mathcal{A}v_2,[a,b,\mathcal{A}\alpha_\mathfrak{h}^{-1}(v_3)]_\mathfrak{g}]_\mathfrak{g}-\mathcal{A}(\theta(\mathcal{A}v_1, [a,b,\mathcal{A}\alpha_\mathfrak{h}^{-1}(v_3)]_\mathfrak{g})v_2-\theta(\mathcal{A}v_2, [a,b,\mathcal{A}\alpha_\mathfrak{h}^{-1}(v_3)]_\mathfrak{g})v_1)\\
&-\mathcal{A}(D(a,b)\alpha_\mathfrak{h}^{-1}(D(\mathcal{A}v_1,\mathcal{A}v_2)v_3+\theta(\mathcal{A}v_2,\mathcal{A}v_3)v_1-\theta(\mathcal{A}v_1,\mathcal{A}v_3)v_2+\kappa [v_1, v_2, v_3]_\mathfrak{h}))\\
&+[a,b,\mathcal{A}\alpha_\mathfrak{h}^{-1}(D(\mathcal{A}v_1,\mathcal{A}v_2)v_3+\theta(\mathcal{A}v_2,\mathcal{A}v_3)v_1-\theta(\mathcal{A}v_1,\mathcal{A}v_3)v_2+\kappa [v_1, v_2, v_3]_\mathfrak{h})]_\mathfrak{g}\\
=&0.
\end{align*}
Therefore, $\delta \Im(a,b)=0.$
\end{proof}

\begin{remark}
  It should be noted that in Proposition \ref{prop:1-cocycle}, the condition that $\Im(a,b)$ is a 1-cocycle  is that Hom-Lie triple systems $\mathfrak{g}$ and $\mathfrak{h}$ are regular, which is different from \cite{Li}.
\end{remark}

\begin{defn}
Let  $\mathcal{A}:\mathfrak{h}\rightarrow \mathfrak{g}$  be  a $\kappa$-weighted  $\mathcal{O}$-operator    from a   Hom-Lie triple system $(\mathfrak{h}, [-, -, -]_\mathfrak{h},\alpha_\mathfrak{h})$ to another   Hom-Lie triple system $(\mathfrak{g}, [-, -, -]_\mathfrak{g},\alpha_\mathfrak{g})$ with respect to an action $\theta$.
Define the set of $n$-cochains   by
      \begin{equation*}
\mathcal{C}^n_{\mathcal{A}}(\mathfrak{h},\mathfrak{g})= \left\{ \begin{array}{ll}
\mathcal{C}_{\mathrm{HLts}}^{2n-1}(\mathfrak{h},  \mathfrak{g}),  &\mbox{  $n\geq 1;$  }\\
$$\{(a,b)\in \wedge^2\mathfrak{g}~|~\alpha_\mathfrak{g}(a)=a,\alpha_\mathfrak{g}(b)=b\},$$ &\mbox{    $ n=0.$  } \end{array}  \right.
\end{equation*}
Define  $\partial_\mathcal{A}:\mathcal{C}^n_{\mathcal{A}}(\mathfrak{h},\mathfrak{g})\rightarrow\mathcal{C}^{n+1}_{\mathcal{A}}(\mathfrak{h},\mathfrak{g})$ by
     \begin{equation*}
\partial_{\mathcal{A}}= \left\{ \begin{array}{ll}
\delta,  &\mbox{  $n\geq 1;$  }\\
$$\Im,$$ &\mbox{    $ n=0.$  } \end{array}  \right.
\end{equation*}
By $\delta\circ\delta=0$ and Proposition \ref{prop:1-cocycle}, we have $(\oplus_{n=0}^{+\infty}\mathcal{C}^n_{\mathcal{A}}(\mathfrak{h},\mathfrak{g}),\partial)$  is a cochain complex.
Denote the set of $n$-cocycles by $\mathcal{Z}^n_{\mathcal{A}}(\mathfrak{h},\mathfrak{g})$, the set of $n$-coboundaries by $\mathcal{B}^n_{\mathcal{A}}(\mathfrak{h},\mathfrak{g})$, and $n$-th cohomology group by
$$\mathcal{H}^n_{\mathcal{A}}(\mathfrak{h},\mathfrak{g})=\frac{\mathcal{Z}^n_{\mathcal{A}}(\mathfrak{h},\mathfrak{g})}{\mathcal{B}^n_{\mathcal{A}}(\mathfrak{h},\mathfrak{g})}, n\geq 1.$$
\end{defn}

\begin{remark}
 The cohomology theory for $\kappa$-weighted $\mathcal{O}$-operators   on Hom-Lie triple systems enjoys certain functorial properties.
 Let  $\mathcal{A}_1,\mathcal{A}_2:\mathfrak{h}\rightarrow \mathfrak{g}$ be  two  $\kappa$-weighted $\mathcal{O}$-operators   from a Hom-Lie triple system $(\mathfrak{h}, [-, -, -]_\mathfrak{h},\alpha_\mathfrak{h})$ to another Hom-Lie triple system $(\mathfrak{g}, [-, -, -]_\mathfrak{g},\alpha_\mathfrak{g})$ with
respect to an action $\theta$, and   $(\varphi_\mathfrak{h},\varphi_\mathfrak{g})$ be a homomorphism from $\mathcal{A}_1$ to $\mathcal{A}_2$
in which $\varphi_\mathfrak{h}$ is invertible. Define a linear map  $\Phi: \mathcal{C}^n_{\mathcal{A}_1}(\mathfrak{h},\mathfrak{g})\rightarrow \mathcal{C}^n_{\mathcal{A}_2}(\mathfrak{h},\mathfrak{g})$  by
 $$\Phi(f)(v_1,\cdots,v_{2n-1})=\varphi_\mathfrak{g}(f(\varphi_\mathfrak{h}^{-1}(v_1),\cdots,\varphi_\mathfrak{h}^{-1}(v_{2n-1}))),$$
 for any $f\in \mathcal{C}^n_{\mathcal{A}_1}(\mathfrak{h},\mathfrak{g})$ and $v_1,\cdots,v_{2n-1}\in \mathfrak{h}.$  Then it is straightforward
to deduce that $\Phi$ is a cochain map from the cochain complex $(\oplus_{n=1}^{+\infty}\mathcal{C}^n_{\mathcal{A}_1}(\mathfrak{h},\mathfrak{g}),\partial_{\mathcal{A}_1})$  to the cochain complex $(\oplus_{n=1}^{+\infty}\mathcal{C}^n_{\mathcal{A}_2}(\mathfrak{h},\mathfrak{g}),\partial_{\mathcal{A}_2})$. Consequently, it induces a homomorphism  $\Phi^*$ from the cohomology group
 $\mathcal{H}^n_{\mathcal{A}_1}(\mathfrak{h},\mathfrak{g})$ to $\mathcal{H}^n_{\mathcal{A}_2}(\mathfrak{h},\mathfrak{g})$.
\end{remark}

\section{  Deformatons of   weighted  $\mathcal{O}$-operators   on Hom-Lie triple systems }  \label{sec: deformations}
\def\theequation{\arabic{section}.\arabic{equation}}
\setcounter{equation} {0}

In this section, we study  linear deformations and higher order deformations of $\kappa$-weighted $\mathcal{O}$-operators  on Hom-Lie triple systems  via the cohomology theory established in the former
section.

First, we use the cohomology constructed to characterize the linear deformations of $\kappa$-weighted $\mathcal{O}$-operators   on Hom-Lie triple systems.

\begin{defn}
Let  $\mathcal{A}:\mathfrak{h}\rightarrow \mathfrak{g}$  be  a $\kappa$-weighted  $\mathcal{O}$-operator   from a   Hom-Lie triple system $(\mathfrak{h}, [-, -, -]_\mathfrak{h},\alpha_\mathfrak{h})$ to another   Hom-Lie triple system $(\mathfrak{g}, [-, -, -]_\mathfrak{g},\alpha_\mathfrak{g})$ with respect to an action $\theta$.
A linear deformation of $\mathcal{A}$ is  a  $\kappa$-weighted  $\mathcal{O}$-operator   of the form $\mathcal{A}_t=\mathcal{A}+t\mathcal{A}_1$, where $\mathcal{A}_1:\mathfrak{h}\rightarrow \mathfrak{g}$ is a linear map and $t$ is a parameter with $t^2=0$.
\end{defn}

Suppose $ \mathcal{A}+t\mathcal{A}_1$ is a linear deformation of   $\mathcal{A}$, direct deduction shows that $\mathcal{A}_1\in \mathcal{C}_\mathcal{A}^{1}(\mathfrak{h},\mathfrak{g})$ is   a 1-cocycle on  $( \mathfrak{h}, \{ -,-,-\}_\mathcal{A},\alpha_\mathfrak{h})$ with coefficients in   $(\mathfrak{g},  \alpha_\mathfrak{g}; \theta_\mathcal{A})$. So the cohomology class of $\mathcal{A}_1$ defines an element in $ \mathcal{H}_\mathcal{A}^{1}(\mathfrak{h},\mathfrak{g})$. Further, the 1-cocycle $ \mathcal{A}_1$ is called the  infinitesimal  of the linear deformation $ \mathcal{A}_t$ of $\mathcal{A}$.

\begin{defn}
Let  $\mathcal{A}:\mathfrak{h}\rightarrow \mathfrak{g}$  be  a  $\kappa$-weighted  $\mathcal{O}$-operator    from a  regular Hom-Lie triple system $(\mathfrak{h}, [-, -, -]_\mathfrak{h},\alpha_\mathfrak{h})$ to another regular  Hom-Lie triple system $(\mathfrak{g}, [-, -, -]_\mathfrak{g},\alpha_\mathfrak{g})$ with respect to an action $\theta$. Two linear deformations $\mathcal{A}_t=\mathcal{A}+t\mathcal{A}_1$ and $\mathcal{A}'_t=\mathcal{A}+t\mathcal{A}'_1$ are called equivalent if   there exist  two elements $a,b\in  \mathfrak{g}$ such that $\alpha_\mathfrak{g}(a)=a,\alpha_\mathfrak{g}(b)=b$ and  the pair $(Id_\mathfrak{g}+t\alpha^{-1}_\mathfrak{g}(\mathcal{L}(a,b)),Id_\mathfrak{h}+t\alpha^{-1}_\mathfrak{h}(D(a,b)))$ is a homomorphism from $\mathcal{A}_t$
to $\mathcal{A}'_t$.
\end{defn}

Suppose $\mathcal{A}_t$ and $\mathcal{A}'_t$ are  equivalent,   then Eq. \eqref{2.13} yields that
$$(Id_\mathfrak{g}+t\alpha^{-1}_\mathfrak{g}(\mathcal{L}(a,b)))\mathcal{A}_tu=\mathcal{A}'_t(Id_\mathfrak{h}+t\alpha^{-1}_\mathfrak{h}(D(a,b)))u, \forall u\in\mathfrak{h}.$$
which means  that
$$\mathcal{A}_1u-\mathcal{A}'_1u=\mathcal{A}(\alpha_\mathfrak{h}^{-1}(D(a,b)u))-\alpha_\mathfrak{g}^{-1}([a,b,\mathcal{A}u]_\mathfrak{g})=\mathcal{A}(D(a,b)\alpha_\mathfrak{h}^{-1}(u))-[a,b,\mathcal{A}\alpha_\mathfrak{g}^{-1}(u)]_\mathfrak{g}.$$

By Proposition \ref{prop:1-cocycle}, we have  $\mathcal{A}_1-\mathcal{A}'_1=\Im(a,b)=\partial_\mathcal{A}(a,b)$. So their cohomology classes are the same in $ \mathcal{H}_\mathcal{A}^{1}(\mathfrak{h},\mathfrak{g})$.

Conversely, any 1-cocycle  $\mathcal{A}_1$  gives rise to the linear deformation $\mathcal{A}+t\mathcal{A}_1$.  To sum up, we
have the following result.
\begin{prop}
 Let  $\mathcal{A}:\mathfrak{h}\rightarrow \mathfrak{g}$  be  a $\kappa$-weighted  $\mathcal{O}$-operator    from a regular  Hom-Lie triple system $(\mathfrak{h}, [-, -, -]_\mathfrak{h},\alpha_\mathfrak{h})$ to another  regular Hom-Lie triple system $(\mathfrak{g}, [-, -, -]_\mathfrak{g},\alpha_\mathfrak{g})$ with respect to an action $\theta$. Then there is a bijection between the set of all
equivalence classes of linear deformation of   $\mathcal{A}$  and the first cohomology group $ \mathcal{H}_\mathcal{A}^{1}(\mathfrak{h},\mathfrak{g})$.
\end{prop}

Next,  we introduce a special cohomology class associated to an
$n$-order deformation of a $\kappa$-weighted  $\mathcal{O}$-operator, and show that an $n$-order deformation
of  a $\kappa$-weighted  $\mathcal{O}$-operator   is extendable if and only if this cohomology class in the second cohomology group vanishes.

\begin{defn}
Let  $\mathcal{A}:\mathfrak{h}\rightarrow \mathfrak{g}$  be  a $\kappa$-weighted  $\mathcal{O}$-operator   from a   Hom-Lie triple system  $(\mathfrak{h}, [-, -, -]_\mathfrak{h},\alpha_\mathfrak{h})$ to another   Hom-Lie triple system  $(\mathfrak{g}, [-, -, -]_\mathfrak{g},\alpha_\mathfrak{g})$ with respect to an action $\theta$.
If $\mathcal{A}_t= \sum_{i=0}^{n}t^i\mathcal{A}_i$ with $\mathcal{A}_0=\mathcal{A}, \mathcal{A}_i\in \mathrm{Hom}(\mathfrak{h}, \mathfrak{g}), i=1,\cdots,n$, defines a $\mathbb{K}[[t]]/(t^{n+1})$-module map
 from   $\mathfrak{h}[[t]]/(t^{n+1})$ to the Hom-Lie triple system $\mathfrak{g}[[t]]/(t^{n+1})$ satisfying
 \begin{align*}
& \mathcal{A}_t\circ \alpha_\mathfrak{h}=\alpha_\mathfrak{g}\circ \mathcal{A}_t, \\
& [\mathcal{A}_tu,\mathcal{A}_tv,\mathcal{A}_tw]_\mathfrak{g}=\mathcal{A}_t(D(\mathcal{A}_tu,\mathcal{A}_tv)w-\theta(\mathcal{A}_tu,\mathcal{A}_tw)v+\theta(\mathcal{A}_tv,\mathcal{A}_tw)u+\kappa[u,v,w]_\mathfrak{h}),
\end{align*}
for any $u,v,w\in \mathfrak{h}$, we say that $\mathcal{A}_t$ is an $n$-order  deformation of $\mathcal{A}$.
\end{defn}

\begin{defn}
Let  $\mathcal{A}:\mathfrak{h}\rightarrow \mathfrak{g}$  be  a $\kappa$-weighted  $\mathcal{O}$-operator    from a   Hom-Lie triple system $(\mathfrak{h}, [-, -, -]_\mathfrak{h},\alpha_\mathfrak{h})$ to another   Hom-Lie triple system  $(\mathfrak{g}, [-, -, -]_\mathfrak{g},\alpha_\mathfrak{g})$ with respect to an action $\theta$.
Let $\mathcal{A}_t= \sum_{i=0}^{n}t^i\mathcal{A}_i$ be an $n$-order  deformation of $\mathcal{A}$. If there is a $\mathcal{A}_{n+1}\in  \mathcal{C}_\mathcal{A}^{1}(\mathfrak{h},\mathfrak{g})$
such that  $\mathcal{A}'_t= \mathcal{A}_t+t^{n+1}\mathcal{A}_{n+1}$  is an  $ (n+1) $-order deformation of  $\mathcal{A}$, then
we say that $\mathcal{A}_t$ is extendable.
\end{defn}

\begin{prop} \label{prop:Obs}
 Let  $\mathcal{A}:\mathfrak{h}\rightarrow \mathfrak{g}$  be  a  $\kappa$-weighted  $\mathcal{O}$-operator    from a   Hom-Lie triple system $(\mathfrak{h}, [-, -, -]_\mathfrak{h},\alpha_\mathfrak{h})$ to another   Hom-Lie triple system $(\mathfrak{g}, [-, -, -]_\mathfrak{g},\alpha_\mathfrak{g})$ with respect to an action $\theta$.  Let $\mathcal{A}_t= \sum_{i=0}^{n}t^i\mathcal{A}_i$ be an $n$-order  deformation of $\mathcal{A}$.
 Then $\mathcal{A}_t$ is extendable if
and only if the cohomology class $ [Obs^n]\in \mathcal{H}_\mathcal{A}^{2}(\mathfrak{h},\mathfrak{g})$  vanishes, where
 \begin{align*}
Obs^n(u_1,u_2,u_3)=&\sum_{\mbox{\tiny$\begin{array}{c}
  i+j+k=n+1\\
   i,j,k\geq 1\end{array}$}}\big([\mathcal{A}_iu_1,\mathcal{A}_ju_2,\mathcal{A}_ku_3]_\mathfrak{g}-\mathcal{A}_i(D(\mathcal{A}_ju_1,\mathcal{A}_ku_2)u_3\\
&-\theta(\mathcal{A}_ju_1,\mathcal{A}_ku_3)u_2+\theta(\mathcal{A}_ju_2,\mathcal{A}_ku_3)u_1)\big)-\kappa\mathcal{A}_{n+1}[u_1,v_2,u_3]_\mathfrak{h}.
\end{align*}
\end{prop}
\begin{proof}
Let $\mathcal{A}'_t= \mathcal{A}_t+t^{n+1}\mathcal{A}_{n+1}$ be the extension of  $\mathcal{A}_t$, then for all $u_1,u_2,u_3\in \mathfrak{h}$
 \begin{align*}
& [\mathcal{A}'_tu_1,\mathcal{A}'_tu_2,\mathcal{A}'_tu_3]_\mathfrak{g}=\mathcal{A}'_t(D(\mathcal{A}'_tu_1,\mathcal{A}'_tu_2)u_3-\theta(\mathcal{A}'_tu_1,\mathcal{A}'_tu_3)u_2+\theta(\mathcal{A}'_tu_2,\mathcal{A}'_tu_3)u_1+\kappa[u_1,v_2,u_3]_\mathfrak{h}).
\end{align*}
Expanding the equation and comparing the coefficients of  $t^n$ yields that:
 \begin{align*}
&\sum_{\mbox{\tiny$\begin{array}{c}
  i+j+k=n+1\\
   i,j,k\geq 1\end{array}$}}\big([\mathcal{A}_iu_1,\mathcal{A}_ju_2,\mathcal{A}_ku_3]_\mathfrak{g}-\mathcal{A}_i(D(\mathcal{A}_ju_1,\mathcal{A}_ku_2)u_3
-\theta(\mathcal{A}_ju_1,\mathcal{A}_ku_3)u_2+\theta(\mathcal{A}_ju_2,\mathcal{A}_ku_3)u_1)\big)\\
&+[\mathcal{A}_{n+1}u_1,\mathcal{A}u_2,\mathcal{A}u_3]_\mathfrak{g}+[\mathcal{A}u_1,\mathcal{A}_{n+1}u_2,\mathcal{A}u_3]_\mathfrak{g}+[\mathcal{A}u_1,\mathcal{A}u_2,\mathcal{A}_{n+1}u_3]_\mathfrak{g}-\mathcal{A}(D(\mathcal{A}_{n+1}u_1,\mathcal{A}u_2)u_3\\
&+D(\mathcal{A}u_1,\mathcal{A}_{n+1}u_2)u_3-\theta(\mathcal{A}_{n+1}u_1,\mathcal{A}u_3)u_2-\theta(\mathcal{A}u_1,\mathcal{A}_{n+1}u_3)u_2+\theta(\mathcal{A}_{n+1}u_2,\mathcal{A}u_3)u_1\\
&+\theta(\mathcal{A}u_2,\mathcal{A}_{n+1}u_3)u_1)-\kappa\mathcal{A}_{n+1}[u_1,v_2,u_3]_\mathfrak{h}=0,
\end{align*}
which is equivalent to $ Obs^n=\partial_\mathcal{A}\mathcal{A}_{n+1}$. Hence the cohomology class $ [Obs^n]\in \mathcal{H}_\mathcal{A}^{2}(\mathfrak{h},\mathfrak{g})$  vanishes.

Conversely, suppose that the cohomology class $ Obs^n$
 vanishes, then there exists a 1-cochain $\mathcal{A}_{n+1}\in  \mathcal{C}_\mathcal{A}^{1}(\mathfrak{h},\mathfrak{g})$ such that $ Obs^n=\partial_\mathcal{A}\mathcal{A}_{n+1}$. Set
 $\mathcal{A}'_t= \mathcal{A}_t+t^{n+1}\mathcal{A}_{n+1}$.
Then $\mathcal{A}'_t$ satisfies
 \begin{align*}
&\sum_{\mbox{\tiny$\begin{array}{c}
  i+j+k=l\\
   i,j,k\geq 1\end{array}$}}\big([\mathcal{A}_iu_1,\mathcal{A}_ju_2,\mathcal{A}_ku_3]_\mathfrak{g}-\mathcal{A}_i(D(\mathcal{A}_ju_1,\mathcal{A}_ku_2)u_3-\theta(\mathcal{A}_ju_1,\mathcal{A}_ku_3)u_2+\theta(\mathcal{A}_ju_2,\mathcal{A}_ku_3)u_1)\big)\\
   &-\kappa\mathcal{A}_{l}[u_1,v_2,u_3]_\mathfrak{h}=0, 0\leq l\leq n+1,
\end{align*}
which implies that $\mathcal{A}'_t$ is an $(n + 1)$-order deformation of  $\mathcal{A}$. Hence it is an extension of $\mathcal{A}_t.$
\end{proof}

\begin{defn}
Let $\mathcal{A}_t= \sum_{i=0}^{n}t^i\mathcal{A}_i$ be an $n$-order  deformation of $\mathcal{A}$. Then the cohomology class $ [Obs^n]\in \mathcal{H}_\mathcal{A}^{2}(\mathfrak{h},\mathfrak{g})$ defined in Proposition  \ref{prop:Obs} is called the obstruction class of $\mathcal{A}_t$ being extendable.
\end{defn}

\begin{coro}
 Let  $\mathcal{A}:\mathfrak{h}\rightarrow \mathfrak{g}$  be  a $\kappa$-weighted  $\mathcal{O}$-operator   from a   Hom-Lie triple system $(\mathfrak{h}, [-, -, -]_\mathfrak{h},\alpha_\mathfrak{h})$ to another   Hom-Lie triple system $(\mathfrak{g}, [-, -, -]_\mathfrak{g},\alpha_\mathfrak{g})$ with respect to an action $\theta$.   If  $ \mathcal{H}_\mathcal{A}^{2}(\mathfrak{h},\mathfrak{g})=0,$  then
every 1-cocycle $\mathcal{A}_{n+1}$ in $ \mathcal{Z}_\mathcal{A}^{1}(\mathfrak{h},\mathfrak{g})$ is the infinitesimal of some $(n+1)$-order deformation of $\mathcal{A}$.
\end{coro}

\section{ Hom-post-Lie triple systems } \label{sec: Hom-post-Lts}
\def\theequation{\arabic{section}.\arabic{equation}}
\setcounter{equation} {0}

In this section, we introduce the notion of  Hom-post-Lie triple system, which is the underlying
algebraic structure of $\kappa$-weighted $\mathcal{O}$-operators. Moreover, we show that
there exists a Hom-Lie triple system structure on a Hom-post-Lie triple system and that an action
can also be obtained via the given Hom-post-Lie triple system.

\begin{defn}
(i) A  Hom-post-Lie triple system $( \mathfrak{h}, \lfloor-,-,-\rfloor_\mathfrak{h}, \{-,-,-\}_\mathfrak{h},\alpha_\mathfrak{h})$ consists of a  Hom-Lie triple system $( \mathfrak{h}, \lfloor-,-,-\rfloor_\mathfrak{h}, \alpha_\mathfrak{h})$ and a
trilinear product $\{-,-,-\}_\mathfrak{h}:  \mathfrak{h}\otimes \mathfrak{h}\otimes \mathfrak{h}\rightarrow  \mathfrak{h}$,  satisfying $\alpha_\mathfrak{h}(\{u,v,w\}_\mathfrak{h})=\{\alpha_\mathfrak{h}(u),\alpha_\mathfrak{h}(v),\alpha_\mathfrak{h}(w)\}_\mathfrak{h}$  such that
\begin{align}
\{\alpha_\mathfrak{h}(u),\alpha_\mathfrak{h}(v),\langle w,s,t\rangle_C\}_\mathfrak{h}=&\{\{u,v,w\}_\mathfrak{h},\alpha_\mathfrak{h}(s),\alpha_\mathfrak{h}(t)\}_\mathfrak{h}-
\{\{u,v,s\}_\mathfrak{h},\alpha_\mathfrak{h}(w),\alpha_\mathfrak{h}(t)\}_\mathfrak{h}\nonumber\\
&+\{\alpha_\mathfrak{h}(w),\alpha_\mathfrak{h}(s),\{u,v,t\}_\mathfrak{h}\}_D,\label{3.1}\\
\{\alpha_\mathfrak{h}(u),\alpha_\mathfrak{h}(v),\{ w,s,t\}_\mathfrak{h}\}_D=&\{\{u,v,w\}_D,\alpha_\mathfrak{h}(s),\alpha_\mathfrak{h}(t)\}_\mathfrak{h}+
\{\alpha_\mathfrak{h}(w),\langle u,v,s\rangle_C,\alpha_\mathfrak{h}(t)\}_\mathfrak{h}\nonumber\\
&+\{\alpha_\mathfrak{h}(w),\alpha_\mathfrak{h}(s),\langle u,v,t\rangle_C\}_\mathfrak{h},\label{3.2}\\
\{\lfloor w,s,t\rfloor_\mathfrak{h}, \alpha_\mathfrak{h}(u),\alpha_\mathfrak{h}(v)\}_\mathfrak{h}=&\lfloor  \alpha_\mathfrak{h}(s),\alpha_\mathfrak{h}(t), \{ w,u,v\}_\mathfrak{h}\rfloor_\mathfrak{h}=0, \label{3.3}
\end{align}
where $u,v,w,s,t\in  \mathfrak{h}$,  operations $\langle -,-,-\rangle_C$ and $\{-,-,-\}_D$ are defined to be
\begin{align}
\langle u,v,w\rangle_C=&\{u,v,w\}_D+\{u,v,w\}_\mathfrak{h}-\{v,u,w\}_\mathfrak{h}+ \lfloor u,v,w\rfloor_\mathfrak{h}, \label{3.4}\\
\{ u,v,w\}_D=&\{w,v,u\}_\mathfrak{h}-\{w,u,v\}_\mathfrak{h}. \label{3.5}
\end{align}

(ii) A homomorphism between two  Hom-post-Lie triple systems  $( \mathfrak{h}_1, \lfloor-,-,-\rfloor_{\mathfrak{h}_1}, \{-,-,-\}_{\mathfrak{h}_1},$ $\alpha_{\mathfrak{h}_1})$ and $( \mathfrak{h}_2, \lfloor-,-,-\rfloor_{\mathfrak{h}_2}, \{-,-,-\}_{\mathfrak{h}_2},\alpha_{\mathfrak{h}_2})$ is a linear map $\varphi: \mathfrak{h}_1\rightarrow \mathfrak{h}_2$ satisfying
$\varphi(\alpha_{\mathfrak{h}_1}(x))=\alpha_{\mathfrak{h}_2}(\varphi(x)), ~~\varphi(\lfloor x, y, z\rfloor_{\mathfrak{h}_1})=\lfloor\varphi(x), \varphi(y),\varphi(z)\rfloor_{\mathfrak{h}_2},~\varphi(\{ x, y, z\}_{\mathfrak{h}_1})=\{\varphi(x), \varphi(y),\varphi(z)\}_{\mathfrak{h}_2},$
for all $x,y,z\in \mathfrak{h}_1.$
\end{defn}

\begin{remark}
 (i) Let $( \mathfrak{h}, \lfloor-,-,-\rfloor_\mathfrak{h}, \{-,-,-\}_\mathfrak{h},\alpha_\mathfrak{h})$  be a Hom-post-Lie triple system.
 On the one hand, if  the   bracket $ \{-,-,-\}_\mathfrak{h}=0,$ we
get that  $( \mathfrak{h}, \lfloor-,-,-\rfloor_\mathfrak{h},  \alpha_\mathfrak{h})$   is a Hom-Lie triple system.
On the other hand,  if the   bracket $\lfloor-,-,$ $-\rfloor_\mathfrak{h}=0$, then  $( \mathfrak{h}, \{-,-,-\}_\mathfrak{h},\alpha_\mathfrak{h})$
becomes a Hom-pre-Lie triple system. Thus,  Hom-post-Lie triple systems are generalizations of both Hom-Lie triple systems and
Hom-pre-Lie triple systems. See~\cite{Li} for more details about Hom-pre-Lie triple systems.

(ii) A post-Lie triple system   is a Hom-post-Lie triple system  with $ \alpha_\mathfrak{h}= \mathrm{id}_\mathfrak{h}$.
\end{remark}

\begin{exam}
Let $(\mathfrak{h}, [-, -, -]_\mathfrak{h},\alpha_\mathfrak{h})$  be a 4-dimensional Hom-Lie triple system with a basis  $\{\varepsilon_1,\varepsilon_2,\varepsilon_3,\varepsilon_4\}$   given by Example \ref{exam:4-dimensional Hom-Lts}.
We define a trilinear product $\{-,-,-\}_\mathfrak{h}:  \mathfrak{h}\otimes \mathfrak{h}\otimes \mathfrak{h}\rightarrow  \mathfrak{h}$ by
$$\{\varepsilon_4,\varepsilon_2,\varepsilon_3\}_\mathfrak{h}=\varepsilon_2,\{\varepsilon_4,\varepsilon_2,\varepsilon_2\}_\mathfrak{h}=\varepsilon_1,\{\varepsilon_3,\varepsilon_2,\varepsilon_4\}_\mathfrak{h}=\{\varepsilon_4,\varepsilon_4,\varepsilon_4\}_\mathfrak{h}=-\varepsilon_2.$$
Then, $( \mathfrak{h}, \lfloor-,-,-\rfloor_\mathfrak{h}, \{-,-,-\}_\mathfrak{h},\alpha_\mathfrak{h})$  is a  Hom-post-Lie triple system.
\end{exam}

\begin{prop} \label{prop:subadjacent Hom-Lts}
 Let $( \mathfrak{h}, \lfloor-,-,-\rfloor_\mathfrak{h}, \{-,-,-\}_\mathfrak{h},\alpha_\mathfrak{h})$ be a   Hom-post-Lie triple system. Then

 (i) the triple $( \mathfrak{h}, \langle -,-,-\rangle_C ,\alpha_\mathfrak{h})$ is a   Hom-Lie triple system,  which is called the adjacent Hom-Lie triple system.

 (ii) the map $\mathfrak{R}$ is an action of the  adjacent  Hom-Lie triple system  $( \mathfrak{h}, \langle -,-,-\rangle_C ,\alpha_\mathfrak{h})$  on
the Hom-Lie triple system $( \mathfrak{h}, \lfloor-,-,-\rfloor_\mathfrak{h},  \alpha_\mathfrak{h})$, where
$$\mathfrak{R}:\mathfrak{h}\otimes\mathfrak{h}\rightarrow \mathrm{End}(\mathfrak{h}), ~(u,v)\mapsto(w\mapsto\{w,u,v\}_\mathfrak{h}),~ \forall u,v,w\in \mathfrak{h}.$$

(iii) the identity map $id:\mathfrak{h}\rightarrow \mathfrak{h}$ is a $\kappa$-weighted $\mathcal{O}$-operator    from  $( \mathfrak{h}, \lfloor-,-,-\rfloor_\mathfrak{h},  \alpha_\mathfrak{h})$ to $( \mathfrak{h}, \langle -,-,-\rangle_C ,\alpha_\mathfrak{h})$ with respect to the action $(\mathfrak{h},\alpha_\mathfrak{h};\mathfrak{R}^\dag)$.
\end{prop}
\begin{proof}
 (i) Obviously, for any $u,v,w\in \mathfrak{h}$, by Eqs. \eqref{2.2}, \eqref{3.4} and \eqref{3.5}, we have  $\langle u,v,w\rangle_C+\langle v,u,w\rangle_C=0$ and $\langle u,v,w\rangle_C+\langle w,u,v\rangle_C+\langle v,w,u\rangle_C=0$. Furthermore, for any $u,v,w,s,t\in \mathfrak{h}$, by Eqs. \eqref{2.3} and \eqref{3.1}-\eqref{3.5}, we get
 \begin{align*}
&\langle\langle s, t, u\rangle_C, \alpha_\mathfrak{h}(v), \alpha_\mathfrak{h}(w)\rangle_C+ \langle\alpha_\mathfrak{h}(u),  \langle s, t, v\rangle_C, \alpha_\mathfrak{h}(w)\rangle_C+ \langle\alpha_\mathfrak{h}(u),  \alpha_\mathfrak{h}(v), \langle s, t, w\rangle_C\rangle_C\\
&-\langle\alpha_\mathfrak{h}(s), \alpha_\mathfrak{h}(t), \langle u, v, w\rangle_C\rangle_C\\
=&\{\alpha_\mathfrak{h}(w), \alpha_\mathfrak{h}(v),\langle s, t, u\rangle_C\}_\mathfrak{h}-\{\alpha_\mathfrak{h}(w), \langle s, t, u\rangle_C, \alpha_\mathfrak{h}(v)\}_\mathfrak{h}
+\{\langle s, t, u\rangle_C, \alpha_\mathfrak{h}(v), \alpha_\mathfrak{h}(w)\}_\mathfrak{h}\\
&-\{\alpha_\mathfrak{h}(v), \langle s, t, u\rangle_C, \alpha_\mathfrak{h}(w)\}_\mathfrak{h}+\lfloor\langle s, t, u\rangle_C, \alpha_\mathfrak{h}(v), \alpha_\mathfrak{h}(w)\rfloor_\mathfrak{h}+ \{ \alpha_\mathfrak{h}(w), \langle s, t, v\rangle_C,\alpha_\mathfrak{h}(u)\}_\mathfrak{h}\\
&- \{ \alpha_\mathfrak{h}(w), \alpha_\mathfrak{h}(u),  \langle s, t, v\rangle_C\}_\mathfrak{h}+ \{\alpha_\mathfrak{h}(u),  \langle s, t, v\rangle_C, \alpha_\mathfrak{h}(w)\}_\mathfrak{h}- \{\langle s, t, v\rangle_C, \alpha_\mathfrak{h}(u),  \alpha_\mathfrak{h}(w)\}_\mathfrak{h}\\
&+\lfloor\alpha_\mathfrak{h}(u),  \langle s, t, v\rangle_C, \alpha_\mathfrak{h}(w)\rfloor_\mathfrak{h}+ \{\langle s, t, w\rangle_C,  \alpha_\mathfrak{h}(v),\alpha_\mathfrak{h}(u)\}_\mathfrak{h}- \{\langle s, t, w\rangle_C, \alpha_\mathfrak{h}(u),  \alpha_\mathfrak{h}(v)\}_\mathfrak{h}\\
&+  \{\alpha_\mathfrak{h}(u),  \alpha_\mathfrak{h}(v), \langle s, t, w\rangle_C\}_\mathfrak{h}- \{\alpha_\mathfrak{h}(v),\alpha_\mathfrak{h}(u),  \langle s, t, w\rangle_C\}_\mathfrak{h}+ \lfloor\alpha_\mathfrak{h}(u),  \alpha_\mathfrak{h}(v), \langle s, t, w\rangle_C\rfloor_\mathfrak{h}\\
&-\{\langle u, v, w\rangle_C, \alpha_\mathfrak{h}(t),\alpha_\mathfrak{h}(s)\}_\mathfrak{h}+\{\langle u, v, w\rangle_C, \alpha_\mathfrak{h}(s), \alpha_\mathfrak{h}(t)\}_\mathfrak{h}-\{\alpha_\mathfrak{h}(s), \alpha_\mathfrak{h}(t), \langle u, v, w\rangle_C\}_\mathfrak{h}\\
&+\{\alpha_\mathfrak{h}(t), \alpha_\mathfrak{h}(s), \langle u, v, w\rangle_C\}_\mathfrak{h}-\lfloor\alpha_\mathfrak{h}(s), \alpha_\mathfrak{h}(t), \langle u, v, w\rangle_C\rfloor_\mathfrak{h}\\
=&0.
\end{align*}
Therefore,  $( \mathfrak{h}, \langle -,-,-\rangle_C ,\alpha_\mathfrak{h})$ is a   Hom-Lie triple system.

(ii) For all $u,v,w\in \mathfrak{h}$, we have
 \begin{align*}
&\mathfrak{L}(u,v)w=\mathfrak{R}(v,u)w-\mathfrak{R}(u,v)w=\{w,v,u\}_\mathfrak{h}-\{w,u,v\}_\mathfrak{h}=\{u,v, w\}_D.
\end{align*}
Obviously, $\mathfrak{R}(\alpha_\mathfrak{h}(u),\alpha_\mathfrak{h}(v))\alpha_\mathfrak{h}(w)=\alpha_\mathfrak{h}(\mathfrak{R}(u,v)(w)).$
Furthermore, for any $u,v,w,s,t\in \mathfrak{h}$, by Eqs. \eqref{2.3} and \eqref{3.1}-\eqref{3.5}, we get
 \begin{align*}
&\mathfrak{R}(\alpha_\mathfrak{h}(u),\alpha_\mathfrak{h}(v))\mathfrak{R}(s,t)w-\mathfrak{R}(\alpha_\mathfrak{h}(t),\alpha_\mathfrak{h}(v))\mathfrak{R}(s,u)w-\mathfrak{R}(\alpha_\mathfrak{h}(s), \langle t,u,v\rangle_C)\alpha_\mathfrak{h}(w)\\
&+\mathfrak{L}(\alpha_\mathfrak{h}(t),\alpha_\mathfrak{h}(u))\mathfrak{R}(s,v)w\\
=&\{\{w,s,t\}_\mathfrak{h},\alpha_\mathfrak{h}(u),\alpha_\mathfrak{h}(v)\}_\mathfrak{h}-\{\{w,s,u\}_\mathfrak{h},\alpha_\mathfrak{h}(t),\alpha_\mathfrak{h}(v)\}_\mathfrak{H}-\{\alpha_\mathfrak{h}(w),\alpha_\mathfrak{h}(s), \langle t,u,v\rangle_C\}_\mathfrak{h}\\
&+\{\alpha_\mathfrak{h}(t),\alpha_\mathfrak{h}(u),\{w,s,v\}_\mathfrak{h}\}_D\\
=&0,\\
& \mathfrak{R}(\alpha_\mathfrak{h}(u),\alpha_\mathfrak{h}(v))\mathfrak{L}(s,t)w-\mathfrak{L}(\alpha_\mathfrak{h}(s),\alpha_\mathfrak{h}(t))\mathfrak{R}(u,v)w+\mathfrak{R}(\langle s,t,u\rangle_C,\alpha_\mathfrak{h}(v))\alpha_\mathfrak{h}(w)\\
&+\mathfrak{R}(\alpha_\mathfrak{h}(u),\langle s,t,v\rangle_C)\alpha_\mathfrak{h}(w)\\
=&\{\{s,t,w\}_D,\alpha_\mathfrak{h}(u),\alpha_\mathfrak{h}(v)\}_\mathfrak{h}-\{\alpha_\mathfrak{h}(s),\alpha_\mathfrak{h}(t),\{w,u,v\}_\mathfrak{h}\}_D+\{\alpha_\mathfrak{h}(w),\langle s,t,u\rangle_C,\alpha_\mathfrak{h}(v)\}_\mathfrak{h}\\
&+\{\alpha_\mathfrak{h}(w),\alpha_\mathfrak{h}(u),\langle s,t,v\rangle_C\}_\mathfrak{h}\\
=&0.
 \end{align*}
Hence,  $(\mathfrak{h}, \alpha_\mathfrak{h}; \mathfrak{R})$ is a representation of the adjacent Hom-Lie triple system  $( \mathfrak{h}, \langle -,-,-\rangle_C ,\alpha_\mathfrak{h})$.
On the other hand,  by Eq.   \eqref{3.3}, we have $  \mathfrak{R}(\alpha_\mathfrak{h}(u),\alpha_\mathfrak{h}(v))\lfloor w,s,t \rfloor_\mathfrak{h}=0,$
$\lfloor\alpha_\mathfrak{h}(s),\alpha_\mathfrak{h}(t), \mathfrak{R}(u,v)w \rfloor_\mathfrak{h}=0.$
So, $\mathfrak{R}$ is an action of the adjacent  Hom-Lie triple system  $( \mathfrak{h}, \langle -,-,-\rangle_C ,\alpha_\mathfrak{h})$  on
the Hom-Lie triple system $( \mathfrak{h}, \lfloor-,-,-\rfloor_\mathfrak{h},  \alpha_\mathfrak{h})$.

(iii) For any  $u,v,w\in \mathfrak{h}$,  by  Eq.   \eqref{3.4}, we have
 \begin{align*}
&id( \mathfrak{L}(id(u),id(v))w- \mathfrak{R}(id(u),id(w))v+\mathfrak{R}(id(v), id(w))u+\kappa\lfloor u,v,w\rfloor_\mathfrak{h})\\
=&\{u,v,w\}_D- \{v,u,w\}_\mathfrak{h}+\{u, v, w\}_\mathfrak{h}+\kappa\lfloor u,v,w\rfloor_\mathfrak{h}\\
=&\langle id(u),id(v),id(w)\rangle_C.
 \end{align*}
Thus,  the identity map $id:\mathfrak{h}\rightarrow \mathfrak{h}$ is a $\kappa$-weighted $\mathcal{O}$-operator.
\end{proof}

\begin{coro}
  Let $\varphi: ( \mathfrak{h}_1, \lfloor-,-,-\rfloor_{\mathfrak{h}_1}, \{-,-,-\}_{\mathfrak{h}_1},\alpha_{\mathfrak{h}_1})\rightarrow( \mathfrak{h}_2, \lfloor-,-,-\rfloor_{\mathfrak{h}_2}, \{-,-,-\}_{\mathfrak{h}_2},\alpha_{\mathfrak{h}_2})$ be a Hom-post-Lie triple system homomorphism. Then, $\varphi$ is
also a Hom-Lie triple system homomorphism between the  subadjacent Hom-Lie triple system
from $( \mathfrak{h}_1, \langle -,-,-\rangle_{C_1} ,\alpha_{\mathfrak{h}_1})$
to $( \mathfrak{h}_2, \langle -,-,-\rangle_{C_2} ,\alpha_{\mathfrak{h}_2})$.

\end{coro}

\begin{prop}
 Let  $\mathcal{A}:\mathfrak{h}\rightarrow \mathfrak{g}$ be  a $\kappa$-weighted $\mathcal{O}$-operator   from a Hom-Lie triple system $(\mathfrak{h}, [-, -, -]_\mathfrak{h},\alpha_\mathfrak{h})$ to another Hom-Lie triple system $(\mathfrak{g}, [-, -, -]_\mathfrak{g},\alpha_\mathfrak{g})$ with
respect to an action $\theta$. Then,

(i) the 4-tuple $(\mathfrak{h}, \lfloor-, -, -\rfloor_\mathfrak{h}, \{-, -, -\}_\mathfrak{h},$ $ \alpha_\mathfrak{h})$ is a  Hom-post-Lie triple system, where
$$\lfloor u, v, w\rfloor_\mathfrak{h}=\kappa[u, v, w]_\mathfrak{h},~\{u,v,w\}_\mathfrak{h}=\theta(\mathcal{A}v,\mathcal{A}w)u, ~\forall u,v,w\in \mathfrak{h}.$$

(ii) the triple $( \mathfrak{h}, \{ -,-,-\}_\mathcal{A},\alpha_\mathfrak{h})$ is a   Hom-Lie triple system,  which is called the descent Hom-Lie triple system, where
 \begin{align}
\{u,v,w\}_\mathcal{A}=D(\mathcal{A}u,\mathcal{A}v)w+\theta(\mathcal{A}v,\mathcal{A}w)u-\theta(\mathcal{A}u,\mathcal{A}w)v+\kappa [u, v, w]_\mathfrak{h}, ~\forall u,v,w\in \mathfrak{h}.\label{3.6}
 \end{align}
Furthermore,  $\mathcal{A}$ is a Hom-Lie triple system homomorphism from the descent Hom-Lie triple system $( \mathfrak{h}, \{ -,-,-\}_\mathcal{A},\alpha_\mathfrak{h})$ to the Hom-Lie triple system $(\mathfrak{g}, [-, -, -]_\mathfrak{g},\alpha_\mathfrak{g})$.
\end{prop}

\begin{proof}
(i) For any $u,v, w, s, t \in \mathfrak{h}$,  first, by Eqs. \eqref{3.3} and \eqref{3.4}, we have
\begin{align}
\{u,v,w\}_D=&\{w,v,u\}_\mathfrak{h}-\{w,u,v\}_\mathfrak{h}=\theta(\mathcal{A}v,\mathcal{A}u)w-\theta(\mathcal{A}u,\mathcal{A}v)w= D(\mathcal{A}u,\mathcal{A}v)w,\nonumber\\
\langle u,v,w\rangle_C=&\{u,v,w\}_D+\{u,v,w\}_\mathfrak{h}-\{v,u,w\}_\mathfrak{h}+\lfloor u,v,w\rfloor_\mathfrak{h}\nonumber\\
=&D(\mathcal{A}u,\mathcal{A}v)w+\theta(\mathcal{A}v,\mathcal{A}w)u-\theta(\mathcal{A}u,\mathcal{A}w)v+\kappa [u, v, w]_\mathfrak{h}. \label{3.7}
\end{align}
Furthermore, by Eqs. \eqref{2.5},\eqref{2.6}, \eqref{2.11}-\eqref{2.14},  we get
\begin{align*}
&\{\{u,v,w\}_\mathfrak{h},\alpha_\mathfrak{h}(s),\alpha_\mathfrak{h}(t)\}_\mathfrak{h}-\{\{u,v,s\}_\mathfrak{h},\alpha_\mathfrak{h}(w),\alpha_\mathfrak{h}(t)\}_\mathfrak{h} +\{\alpha_\mathfrak{h}(w),\alpha_\mathfrak{h}(s),\{u,v,t\}_\mathfrak{h}\}_D\\
&-\{\alpha_\mathfrak{h}(u),\alpha_\mathfrak{h}(v),\langle w,s,t\rangle_C\}_\mathfrak{h}\\
=&\theta(\mathcal{A}\alpha_\mathfrak{h}(s),\mathcal{A}\alpha_\mathfrak{h}(t))\theta(\mathcal{A}v,\mathcal{A}w)u-\theta(\mathcal{A}\alpha_\mathfrak{h}(w),\mathcal{A}\alpha_\mathfrak{h}(t))\theta(\mathcal{A}v,\mathcal{A}s)u +D(\mathcal{A}\alpha_\mathfrak{h}(w),\mathcal{A}\alpha_\mathfrak{h}(s))\theta(\mathcal{A}v,\mathcal{A}t)u\\
&-\theta(\mathcal{A}\alpha_\mathfrak{h}(v),\mathcal{A}(D(\mathcal{A}w,\mathcal{A}s)t+\theta(\mathcal{A}s,\mathcal{A}t)w-\theta(\mathcal{A}w,\mathcal{A}t)s+\kappa [w, s, t]_\mathfrak{h}))\alpha_\mathfrak{h}(u)\\
=&\theta(\alpha_\mathfrak{h}(\mathcal{A}s),\alpha_\mathfrak{h}(\mathcal{A}t))\theta(\mathcal{A}v,\mathcal{A}w)u-\theta(\alpha_\mathfrak{h}(\mathcal{A}w),\alpha_\mathfrak{h}(\mathcal{A}t))\theta(\mathcal{A}v,\mathcal{A}s)u +D(\alpha_\mathfrak{h}(\mathcal{A}w),\alpha_\mathfrak{h}(\mathcal{A}s))\theta(\mathcal{A}v,\mathcal{A}t)u\\
&-\theta(\alpha_\mathfrak{h}(\mathcal{A}v),[\mathcal{A}w, \mathcal{A}s, \mathcal{A}t]_\mathfrak{h})\alpha_\mathfrak{h}(u)\\
=&0,\\
&\{\{u,v,w\}_D,\alpha_\mathfrak{h}(s),\alpha_\mathfrak{h}(t)\}_\mathfrak{h}+\{\alpha_\mathfrak{h}(w),\langle u,v,s\rangle_C,\alpha_\mathfrak{h}(t)\}_\mathfrak{h}+\{\alpha_\mathfrak{h}(w),\alpha_\mathfrak{h}(s),\langle u,v,t\rangle_C\}_\mathfrak{h}\\
&-\{\alpha_\mathfrak{h}(u),\alpha_\mathfrak{h}(v),\{ w,s,t\}_\mathfrak{h}\}_D\\
=&\theta(\alpha_\mathfrak{h}(\mathcal{A}s),\alpha_\mathfrak{h}(\mathcal{A}t))D(\mathcal{A}u,\mathcal{A}v)w+\theta([\mathcal{A}u,\mathcal{A}v,\mathcal{A}s]_\mathfrak{h},\alpha_\mathfrak{h}(\mathcal{A}t))\alpha_\mathfrak{h}(w)
+\theta(\alpha_\mathfrak{h}(\mathcal{A}s),[ \mathcal{A}u,\mathcal{A}v,\mathcal{A}t]_\mathfrak{h})\alpha_\mathfrak{h}(w)\\
&-D(\alpha_\mathfrak{h}(\mathcal{A}u),\alpha_\mathfrak{h}(\mathcal{A}v))\theta(\mathcal{A}s,\mathcal{A}t)w\\
=&0,\\
&\{\lfloor w, s, t\rfloor_\mathfrak{h}, \alpha_\mathfrak{h}(u), \alpha_\mathfrak{h}(v)\}_\mathfrak{h}=\kappa\theta(\alpha_\mathfrak{g}(\mathcal{A}u),\alpha_\mathfrak{g}(\mathcal{A}v))[w, s, t]_\mathfrak{h}\\
=&0,\\
&\lfloor\alpha_\mathfrak{h}(u), \alpha_\mathfrak{h}(v), \{w, s, t\}_\mathfrak{h}\rfloor_\mathfrak{h}=\kappa[\alpha_\mathfrak{h}(u), \alpha_\mathfrak{h}(v), \theta(\mathcal{A}s,\mathcal{A}t)w]_\mathfrak{h}\\
=&0.
\end{align*}
So  $(\mathfrak{h}, \lfloor-, -, -\rfloor_\mathfrak{h}, \{-, -, -\}_\mathfrak{h},$ $ \alpha_\mathfrak{h})$ is a  Hom-post-Lie triple system.

(ii) From  Eq. \eqref{3.7}  and Proposition \ref{prop:subadjacent Hom-Lts},  we have $( \mathfrak{h}, \{ -,-,-\}_\mathcal{A},\alpha_\mathfrak{h})$ is a   Hom-Lie triple system. Furthermore, by Eqs. \eqref{2.11} and \eqref{2.12}, $\mathcal{A}$ is a Hom-Lie triple system homomorphism from   $( \mathfrak{h}, \{ -,-,-\}_\mathcal{A},\alpha_\mathfrak{h})$ to   $(\mathfrak{g}, [-, -, -]_\mathfrak{g},\alpha_\mathfrak{g})$.
\end{proof}
\begin{prop}
  Let  $\mathcal{A}_1,\mathcal{A}_2:\mathfrak{h}\rightarrow \mathfrak{g}$ be  two  $\kappa$-weighted $\mathcal{O}$-operators   from a Hom-Lie triple system $(\mathfrak{h}, [-, -, -]_\mathfrak{h},\alpha_\mathfrak{h})$ to another Hom-Lie triple system $(\mathfrak{g}, [-, -, -]_\mathfrak{g},\alpha_\mathfrak{g})$ with
respect to an action $\theta$, and   $(\varphi_\mathfrak{h},\varphi_\mathfrak{g})$ a homomorphism from $\mathcal{A}_1$ to $\mathcal{A}_2$.
Let $(\mathfrak{h}, \lfloor-, -, -\rfloor_{\mathfrak{h}}, \{-, -, $ $-\}_{\mathfrak{h}_1},  \alpha_\mathfrak{h})$ and $(\mathfrak{h}, \lfloor-, -, -\rfloor_{\mathfrak{h}},$ $ \{-, -, -\}_{\mathfrak{h}_2},  \alpha_\mathfrak{h})$ be the induced Hom-post-Lie triple systems.
Then,  $\varphi_\mathfrak{h}$ is a homomorphism from the Hom-post-Lie triple system $(\mathfrak{h}, \lfloor-, -, -\rfloor_{\mathfrak{h}},$ $ \{-, -, -\}_{\mathfrak{h}_1},  \alpha_\mathfrak{h})$ to   $(\mathfrak{h}, \lfloor-, -, -\rfloor_{\mathfrak{h}}, \{-, -, -\}_{\mathfrak{h}_2},  \alpha_\mathfrak{h})$.
\end{prop}

\begin{proof}
For any $u,v,w\in \mathfrak{h}$, by Eqs. \eqref{2.13} and \eqref{2.14}, we have
\begin{align*}
 \varphi_\mathfrak{h}(\{u,v,w\}_{\mathfrak{h}_1})=&\varphi_\mathfrak{h}(\theta(\mathcal{A}_1v,\mathcal{A}_1w)u)=\theta(\varphi_\mathfrak{g}(\mathcal{A}_1v),\varphi_\mathfrak{g}(\mathcal{A}_1w))\varphi_\mathfrak{h}(u)\\
 =&\theta(\mathcal{A}_2\varphi_\mathfrak{h}(v),\mathcal{A}_2\varphi_\mathfrak{h}(w))\varphi_\mathfrak{h}(u)\\
 =&\{\varphi_\mathfrak{h}(u), \varphi_\mathfrak{h}(v), \varphi_\mathfrak{h}(w)\}_{\mathfrak{h}_2}.
\end{align*}
And because  $\varphi_\mathfrak{h}$ is a Hom-Lie triple system homomorphism on $(\mathfrak{h}, [-, -, -]_\mathfrak{h},\alpha_\mathfrak{h})$. Hence,  $\varphi_\mathfrak{h}$ is a homomorphism from  $(\mathfrak{h}, \lfloor-, -, -\rfloor_{\mathfrak{h}}, \{-, -, -\}_{\mathfrak{h}_1},  \alpha_\mathfrak{h})$ to  $(\mathfrak{h}, \lfloor-, -, -\rfloor_{\mathfrak{h}}, \{-, -, -\}_{\mathfrak{h}_2},  \alpha_\mathfrak{h})$.
\end{proof}

\section{ From     Hom-post-Lie algebras  to Hom-post-Lie triple systems } \label{sec: Hom-post-Lie}
\def\theequation{\arabic{section}.\arabic{equation}}
\setcounter{equation} {0}

In this section, we establish connections between  $\kappa$-weighted $\mathcal{O}$-operators   on Hom-Lie algebras
and their corresponding Hom-Lie triple systems. Furthermore, we can construct
a Hom-post-Lie triple system starting from a Hom-post-Lie algebra.

First, we recall the representation of Hom-Lie  algebra from \cite{Sheng12}.

 A  representation of a Hom-Lie   algebra $( \mathfrak{g}, [-,-]_\mathfrak{g},  \alpha_\mathfrak{g})$ on a Hom-vector space  $(\mathfrak{h},   \alpha_\mathfrak{h})$  is a   linear map $\rho:  \mathfrak{g}\rightarrow \mathrm{End}( \mathfrak{h})$, for all $x,y\in  \mathfrak{h}$, such that
\begin{align}
&\rho(\alpha_\mathfrak{g}(x))\circ\alpha_\mathfrak{h}=\alpha_\mathfrak{h}\circ\rho(x),  \label{6.1}\\
&\rho([x,y]_\mathfrak{g})\circ\alpha_\mathfrak{h}=\rho(\alpha_\mathfrak{g}(x))\circ\rho(y)-\rho(\alpha_\mathfrak{g}(y)))\circ\rho(x). \label{6.2}
\end{align}

A  action  of a Hom-Lie  algebra $( \mathfrak{g}, [-,-]_\mathfrak{g},  \alpha_\mathfrak{g})$ on another Hom-Lie  algebra  $(\mathfrak{h}, [-,-]_\mathfrak{h},   \alpha_\mathfrak{h})$  is a   linear map $\rho:  \mathfrak{g}\rightarrow \mathrm{End}(\mathfrak{h})$ satisfying Eqs  \eqref{6.1}, \eqref{6.2} and
\begin{align}
\rho(\alpha_\mathfrak{g}(x))[u,v]_\mathfrak{h}=&[\rho(x)u,\alpha_\mathfrak{h}(v)]_\mathfrak{h}+[\alpha_\mathfrak{h}(u),\rho(x)v]_\mathfrak{h},\label{6.3}\\
[\rho(x)u,\alpha_\mathfrak{h}(v)]_\mathfrak{h}=&0, \label{6.4}
\end{align}
for all $x\in  \mathfrak{g}$ and $u,v\in \mathfrak{h}$.

Let $(\mathfrak{h},   \alpha_\mathfrak{h}; \rho)$ be a representation of a Hom-Lie   algebra $( \mathfrak{g}, [-,-]_\mathfrak{g},  \alpha_\mathfrak{g})$, then
$(\mathfrak{h},   \alpha_\mathfrak{h}; \theta_\rho)$ is also a representation of the Hom-Lie triple system $( \mathfrak{g}, [-,-,-]_\mathfrak{g},  \alpha_\mathfrak{g})$,
where  $$\theta_\rho(x,y)=\rho(\alpha_\mathfrak{h}(y))\rho(x), [x,y,z]_{\mathfrak{g}}=[[x,y]_\mathfrak{g},\alpha_\mathfrak{g}(z)]_\mathfrak{g}, \forall x,y,z\in \mathfrak{g}.$$
Furthermore, if $\rho$ is  an action  of a Hom-Lie  algebra $( \mathfrak{g}, [-,-]_\mathfrak{g},  \alpha_\mathfrak{g})$ on another Hom-Lie  algebra  $(\mathfrak{h}, [-,-]_\mathfrak{h},   \alpha_\mathfrak{h})$,
then  $\theta_\rho$ is also an action  of a Hom-Lie triple system $( \mathfrak{g}, [-,-,-]_\mathfrak{g},  \alpha_\mathfrak{g})$ on another Hom-Lie triple system  $(\mathfrak{h}, [-,-,-]_\mathfrak{h},   \alpha_\mathfrak{h})$.

\begin{defn}
A $\kappa$-weighted $\mathcal{O}$-operator     from a Hom-Lie algebra   $( \mathfrak{h}, [-,-]_\mathfrak{h},  \alpha_\mathfrak{h})$    to another   Hom-Lie algebra   $( \mathfrak{g}, [-,-]_\mathfrak{g},  \alpha_\mathfrak{g})$ with respect to an action $\rho$  is a linear map $\mathcal{A}: \mathfrak{h}\rightarrow  \mathfrak{g}$
satisfying the following equations:
\begin{align}
 \mathcal{A}\circ\alpha_\mathfrak{h}=&\alpha_\mathfrak{g}\circ  \mathcal{A},\label{6.5}\\
[\mathcal{A} u, \mathcal{A}v]_\mathfrak{g}=& \mathcal{A}(\rho( \mathcal{A}u)v-\rho( \mathcal{A}v)u+\kappa[\mu,\nu]_\mathfrak{h}), \label{6.6}
\end{align}
for any  $u,u\in  \mathfrak{h}$.
\end{defn}

\begin{remark}
 Let $\mathcal{A}$ be a $\kappa$-weighted $\mathcal{O}$-operator    from a Hom-Lie algebra   $( \mathfrak{h}, [-,-]_\mathfrak{h},  \alpha_\mathfrak{h})$    to another   Hom-Lie algebra   $( \mathfrak{g}, [-,-]_\mathfrak{g},  \alpha_\mathfrak{g})$ with respect to an action $\rho$. Then $\mathcal{A}$ is also a $\kappa$-weighted  $\mathcal{O}$-operator   from a Hom-Lie triple system   $( \mathfrak{h}, [-,-,-]_\mathfrak{h},  \alpha_\mathfrak{h})$    to another   Hom-Lie triple system   $( \mathfrak{g}, [-,-,-]_\mathfrak{g},  \alpha_\mathfrak{g})$ with respect to an action $\theta_\rho$.
\end{remark}

\begin{defn}
 A  Hom-post-Lie algebra $(\mathfrak{h},[-,-],\star, \alpha_\mathfrak{h})$ is a Hom-Lie algebra $(\mathfrak{h},[-,-], \alpha_\mathfrak{h})$
together with a bilinear map $\star: \mathfrak{h}\otimes \mathfrak{h}\rightarrow\mathfrak{h}$, for any $u,v,w\in\mathfrak{h}$, satisfying $\alpha_\mathfrak{h}(u\star v)=\alpha_\mathfrak{h}(u)\star\alpha_\mathfrak{h}(v)$, such that:
\begin{align}
\alpha_\mathfrak{h}(u)\star [v,w]=&[u\star v, \alpha_\mathfrak{h}(w)]+[\alpha_\mathfrak{h}(v),u\star w],\label{6.7}\\
([u,v]+u\star v-v\star u)\star\alpha_\mathfrak{h}(w)=&\alpha_\mathfrak{h}(u)\star (v\star w)-\alpha_\mathfrak{h}(v)\star(u\star w),\label{6.8}\\
\alpha_\mathfrak{h}(u)\star [v,w] =&[u\star v,\alpha_\mathfrak{h}(w)] =0.\label{6.9}
\end{align}
\end{defn}
It is noted that the concept of Hom-post-Lie algebra introduced by us is different from that in \cite{Bakayoko}.
If the   bracket $[-,-]=0$, then  $( \mathfrak{h}, \star,\alpha_\mathfrak{h})$
becomes a Hom-pre-Lie algebra.   See~\cite{Liu,Sun17} for more details about Hom-pre-Lie algebras.

Direct verification, we have the following result.
\begin{prop}
 Let  $(\mathfrak{h},[-,-],\star, \alpha_\mathfrak{h})$  be a Hom-post-Lie algebra. Then $(\mathfrak{h},[-,-]_C, \alpha_\mathfrak{h})$ is a Hom-Lie algebra called the adjacent Hom-Lie algebra,
 where
 \begin{align}
 [u,v]_C=u\star v-v\star u+[u,v].\label{6.10}
 \end{align}
\end{prop}

\begin{prop}
 Let $\mathcal{A}: \mathfrak{h}\rightarrow  \mathfrak{g}$  be a $\kappa$-weighted  $\mathcal{O}$-operator    from a Hom-Lie algebra   $( \mathfrak{h}, [-,-]_\mathfrak{h},  \alpha_\mathfrak{h})$    to another   Hom-Lie algebra   $( \mathfrak{g}, [-,-]_\mathfrak{g},  \alpha_\mathfrak{g})$ with respect to an action $\rho$. Then,

(i) the 4-tuple $(\mathfrak{h}, [-, -], \star, \alpha_\mathfrak{h})$ is a  Hom-post-Lie algebra, where
$$[ u, v]=\kappa[u, v]_\mathfrak{h},~ u\star v =\rho(\mathcal{A}u)v, ~\forall u,v,w\in \mathfrak{h}.$$

(ii) the triple $( \mathfrak{h}, [-,-]_\mathcal{A},\alpha_\mathfrak{h})$ is a   Hom-Lie algebra,  which is called the descent Hom-Lie algebra, where
 \begin{align*}
[u,v]_\mathcal{A}=\rho( \mathcal{A}u)v-\rho( \mathcal{A}v)u+\kappa[\mu,\nu]_\mathfrak{h}, ~\forall u,v\in \mathfrak{h}.
 \end{align*}
Furthermore,  $\mathcal{A}$ is a Hom-Lie algebra homomorphism from the descent Hom-Lie algebra $( \mathfrak{h}, [-,-]_\mathcal{A},\alpha_\mathfrak{h})$ to the Hom-Lie algebra $(\mathfrak{g}, [-, -]_\mathfrak{g},\alpha_\mathfrak{g})$.

(iii) the map $ \mathrm{ad}$ is an action of the descent Hom-Lie algebra $( \mathfrak{h}, [-,-]_\mathcal{A},\alpha_\mathfrak{h})$  on
the   Hom-Lie algebra $( \mathfrak{h}, [-,-]_\mathfrak{h},  \alpha_\mathfrak{h})$, where
$$\mathrm{ad}:\mathfrak{h}\rightarrow \mathrm{End}(\mathfrak{h}), ~u\mapsto(v\mapsto u\star v),~ \forall u,v\in \mathfrak{h}.$$
\end{prop}
\begin{proof}
(i) For any $u,v, w\in \mathfrak{h}$,  by  Eqs  \eqref{6.1}-\eqref{6.6},  we  have
\begin{align*}
&\alpha_\mathfrak{h}(u)\star [v,w]-[u\star v, \alpha_\mathfrak{h}(w)]-[\alpha_\mathfrak{h}(v),u\star w]\\
=& \kappa\rho(\alpha_\mathfrak{g}(\mathcal{A}u))[v, w]_\mathfrak{h}-\kappa[\rho(\mathcal{A}u)v, \alpha_\mathfrak{h}(w)]_\mathfrak{h}-\kappa[\alpha_\mathfrak{h}(v),\rho(\mathcal{A}u)w]_\mathfrak{h}\\
=&0,\\
&([u,v]+u\star v-v\star u)\star\alpha_\mathfrak{h}(w)-\alpha_\mathfrak{h}(u)\star (v\star w)+\alpha_\mathfrak{h}(v)\star(u\star w),\\
=&\rho(\mathcal{A}(\kappa[u, v]_\mathfrak{h}+\rho(\mathcal{A}u)v-\rho(\mathcal{A}v)u)\alpha_\mathfrak{h}(w)-\rho(\alpha_\mathfrak{g}(\mathcal{A}u))\rho(\mathcal{A}v)w+\rho(\alpha_\mathfrak{g}(\mathcal{A}v))\rho(\mathcal{A}u)w\\
=&\rho([\mathcal{A}u,\mathcal{A}v]_\mathfrak{g})\alpha_\mathfrak{h}(w)-\rho(\alpha_\mathfrak{g}(\mathcal{A}u))\rho(\mathcal{A}v)w+\rho(\alpha_\mathfrak{g}(\mathcal{A}v))\rho(\mathcal{A}u)w\\
=&0,\\
&\alpha_\mathfrak{h}(u)\star [v,w]= \kappa\rho(\alpha_\mathfrak{g}(\mathcal{A}u))[v,w]_\mathfrak{h}=0,\\
&[u\star v,\alpha_\mathfrak{h}(w)] =\kappa [\rho(\mathcal{A}u)v,\alpha_\mathfrak{h}(w)]_\mathfrak{h}=0.
\end{align*}
So  $(\mathfrak{h}, \lfloor-, -, -\rfloor_\mathfrak{h}, \{-, -, -\}_\mathfrak{h},$ $ \alpha_\mathfrak{h})$ is a   Hom-post-Lie algebra.

(ii) By direct calculation, we have $( \mathfrak{h}, [-,-]_\mathcal{A},\alpha_\mathfrak{h})$ is a   Hom-Lie algebra and  $\mathcal{A}:( \mathfrak{h}, [-,-]_\mathcal{A},\alpha_\mathfrak{h})\rightarrow(\mathfrak{g}, [-, -]_\mathfrak{g},\alpha_\mathfrak{g})$ is a Hom-Lie algebra homomorphism.

(iii) For any $u,v, w\in \mathfrak{h}$,  it is obvious that $\mathrm{ad}$ is  a  representation of the descent Hom-Lie   algebra $( \mathfrak{h}, [-,-]_\mathcal{A},\alpha_\mathfrak{h})$ on $( \mathfrak{h},    \alpha_\mathfrak{h})$.
Further, we have
\begin{align*}
 \mathrm{ad}(\alpha_\mathfrak{h}(w))[u,v]_\mathfrak{h}=&\alpha_\mathfrak{h}(w)[u,v]_\mathfrak{h}=\rho(\alpha_\mathfrak{g}(\mathcal{A}w))[u,v]_\mathfrak{h}=[\rho(\mathcal{A}w)u,\alpha_\mathfrak{h}(v)]_\mathfrak{h}+[\alpha_\mathfrak{h}(u),\rho(\mathcal{A}w)v]_\mathfrak{h}\\
 =&[ \mathrm{ad}(w)u,\alpha_\mathfrak{h}(v)]_\mathfrak{h}+[\alpha_\mathfrak{h}(u), \mathrm{ad}(w)v]_\mathfrak{h}, \\
[ \mathrm{ad}(w)u,\alpha_\mathfrak{h}(v)]_\mathfrak{h}=&[\rho(\mathcal{A}w)u,\alpha_\mathfrak{h}(v)]_\mathfrak{h}=0.
\end{align*}
Therefore, the map $ \mathrm{ad}$ is an action of  $( \mathfrak{h}, [-,-]_\mathcal{A},\alpha_\mathfrak{h})$  on
$( \mathfrak{h}, [-,-]_\mathfrak{h},  \alpha_\mathfrak{h})$.
\end{proof}

\begin{prop}
 Let  $(\mathfrak{h},[-,-],\star, \alpha_\mathfrak{h})$  is a Hom-post-Lie algebra. We define two  trilinear brackets $\lfloor-,-,-\rfloor_{\mathfrak{h}}, \{-,-,-\}_{\mathfrak{h}}: \mathfrak{h}\otimes \mathfrak{h}\otimes \mathfrak{h}\rightarrow  \mathfrak{h} $ by
 $$\lfloor u,v,w\rfloor_{\mathfrak{h}}=[[u,v],\alpha_\mathfrak{h}(w)],~~\{u,v,w\}_{\mathfrak{h}}= \alpha_\mathfrak{h}(w) \star (v\star u),~~\forall u,v,w\in \mathfrak{h}.$$
 Then, $( \mathfrak{h}, \lfloor-,-,-\rfloor_\mathfrak{h}, \{-,-,-\}_\mathfrak{h},\alpha_\mathfrak{h})$ is a   Hom-post-Lie triple system.
\end{prop}

\begin{proof}
For any $u,v,w,s,t\in  \mathfrak{h}$, First, by  Eqs  \eqref{6.7}-\eqref{6.10}, we have
\begin{align*}
&\{ u,v,w\}_D=\{w,v,u\}_\mathfrak{h}-\{w,u,v\}_\mathfrak{h}=\alpha_\mathfrak{h}(u) \star (v\star w)-\alpha_\mathfrak{h}(v) \star (u\star w),\\
&[[u,v]_C,\alpha_\mathfrak{h}(w)]_C=[u\star v-v\star u+[u,v],\alpha_\mathfrak{h}(w)]_C\\
=&[u\star v,\alpha_\mathfrak{h}(w)]_C-[v\star u,\alpha_\mathfrak{h}(w)]_C+[[u,v],\alpha_\mathfrak{h}(w)]_C\\
=&(u\star v)\star \alpha_\mathfrak{h}(w)-\alpha_\mathfrak{h}(w)\star (u\star v)+[u\star v,\alpha_\mathfrak{h}(w)]-(v\star u)\star \alpha_\mathfrak{h}(w)+\alpha_\mathfrak{h}(w)\star (v\star u)\\
&-[v\star u,\alpha_\mathfrak{h}(w)]+[u,v]\star \alpha_\mathfrak{h}(w)-\alpha_\mathfrak{h}(w)\star [u,v]+[[u,v],\alpha_\mathfrak{h}(w)]\\
=&(u\star v)\star \alpha_\mathfrak{h}(w)-\alpha_\mathfrak{h}(w)\star (u\star v)-(v\star u)\star \alpha_\mathfrak{h}(w)+\alpha_\mathfrak{h}(w)\star (v\star u)\\
&+[u,v]\star \alpha_\mathfrak{h}(w)+[[u,v],\alpha_\mathfrak{h}(w)]\\
=& \alpha_\mathfrak{h}(u)\star(v\star w)-\alpha_\mathfrak{h}(v)\star(u\star w)-\alpha_\mathfrak{h}(w)\star (u\star v) +\alpha_\mathfrak{h}(w)\star (v\star u)+[[u,v],\alpha_\mathfrak{h}(w)]\\
=&\{u,v,w\}_D+\{u,v,w\}_\mathfrak{h}-\{v,u,w\}_\mathfrak{h}+\lfloor u,v,w\rfloor_\mathfrak{h}\\
=&\langle u,v,w\rangle_C.
\end{align*}
Further, we can get
\begin{align*}
&\{\lfloor w,s,t\rfloor_\mathfrak{h}, \alpha_\mathfrak{h}(u),\alpha_\mathfrak{h}(v)\}_\mathfrak{h}=\alpha^2_\mathfrak{h}(v)\star (\alpha_\mathfrak{h}(u)\star [[w,s],\alpha_\mathfrak{h}(t)])  \\
=&0, \\
&\lfloor  \alpha_\mathfrak{h}(s),\alpha_\mathfrak{h}(t), \{ w,u,v\}_\mathfrak{h}\rfloor_\mathfrak{h}=[[  \alpha_\mathfrak{h}(s),\alpha_\mathfrak{h}(t)], \alpha_\mathfrak{h}(v) \star (u\star w)]=-[ \alpha_\mathfrak{h}(v) \star (u\star w),[  \alpha_\mathfrak{h}(s),\alpha_\mathfrak{h}(t)]]\\
=&0,\\
&\{\alpha_\mathfrak{h}(u),\alpha_\mathfrak{h}(v),\langle w,s,t\rangle_C\}_\mathfrak{h}-\{\{u,v,w\}_\mathfrak{h},\alpha_\mathfrak{h}(s),\alpha_\mathfrak{h}(t)\}_\mathfrak{h}+
\{\{u,v,s\}_\mathfrak{h},\alpha_\mathfrak{h}(w),\alpha_\mathfrak{h}(t)\}_\mathfrak{h} \\
&-\{\alpha_\mathfrak{h}(w),\alpha_\mathfrak{h}(s),\{u,v,t\}_\mathfrak{h}\}_D\\
=&\alpha_\mathfrak{h}([[w,s]_C,\alpha_\mathfrak{h}(t)]_C) \star (\alpha_\mathfrak{h}(v)\star\alpha_\mathfrak{h}(u))- \alpha^2_\mathfrak{h}(t)\star(\alpha_\mathfrak{h}(s) \star( \alpha_\mathfrak{h}(w) \star (v\star u)))\\
&+\alpha^2_\mathfrak{h}(t)\star (\alpha_\mathfrak{h}(w)\star(\alpha_\mathfrak{h}(s)\star (v\star u)))- \alpha^2_\mathfrak{h}(w)\star (\alpha_\mathfrak{h}(s) \star(\alpha_\mathfrak{h}(t)\star (v\star u)))\\
&+ \alpha^2_\mathfrak{h}(s)\star (\alpha_\mathfrak{h}(w) \star(\alpha_\mathfrak{h}(t)\star (v\star u)))\\
=&\alpha^2_\mathfrak{h}([w,s]_C)\star (\alpha_\mathfrak{h}(t)\star ( v \star u ))-\alpha^2_\mathfrak{h}(t)\star(\alpha_\mathfrak{h}([w,s]_C)\star ( v \star u ))\\
&- \alpha^2_\mathfrak{h}(t)\star(\alpha_\mathfrak{h}(s) \star( \alpha_\mathfrak{h}(w) \star (v\star u)))+
\alpha^2_\mathfrak{h}(t)\star (\alpha_\mathfrak{h}(w)\star(\alpha_\mathfrak{h}(s)\star (v\star u)))\\
&- \alpha^2_\mathfrak{h}(w)\star (\alpha_\mathfrak{h}(s) \star(\alpha_\mathfrak{h}(t)\star (v\star u)))+ \alpha^2_\mathfrak{h}(s)\star (\alpha_\mathfrak{h}(w) \star(\alpha_\mathfrak{h}(t)\star (v\star u)))\\
=&0.
\end{align*}
Similarly, we also have
\begin{align*}
&\{\alpha_\mathfrak{h}(u),\alpha_\mathfrak{h}(v),\{ w,s,t\}_\mathfrak{h}\}_D-\{\{u,v,w\}_D,\alpha_\mathfrak{h}(s),\alpha_\mathfrak{h}(t)\}_\mathfrak{h}-
\{\alpha_\mathfrak{h}(w),\langle u,v,s\rangle_C,\alpha_\mathfrak{h}(t)\}_\mathfrak{h} \\
&-\{\alpha_\mathfrak{h}(w),\alpha_\mathfrak{h}(s),\langle u,v,t\rangle_C\}_\mathfrak{h}=0.
\end{align*}
Therefore, $( \mathfrak{h}, \lfloor-,-,-\rfloor_\mathfrak{h}, \{-,-,-\}_\mathfrak{h},\alpha_\mathfrak{h})$ is a   Hom-post-Lie triple system.
\end{proof}

We have shown that Hom-Lie algebras,  Hom-post-Lie algebras, Hom-Lie triple systems and  Hom-post-Lie triple systems are
closely related in the sense of commutative diagram of categories as follows:
 $$\aligned
\xymatrix{
  \text{Hom-post-Lie  algebras} \ar[r]  \ar[d] & \text{Hom-post-Lie  triple  systems} \ar[d]\\
  \ar[u] \text{Hom-Lie  algebras and action} \ar[r]& \text{Hom-Lie  triple  systems and action} \ar[u].}
 \endaligned$$

{{\bf Acknowledgments.}  The paper is  supported by the  Foundation of Science and Technology of Guizhou Province(Grant Nos. [2018]1020,  ZK[2022]031, ZK[2023]025),   the National Natural Science Foundation of China (Grant No. 12161013).


\begin{thebibliography}{99}


\bibitem{Ammar} F. Ammar,  Z. Ejbehi,  A.Makhlouf, Cohomology and deformations of Hom-algebras, {\it J. Lie Theory}, {\bf21}(4), (2011), 813--836.


\bibitem{Bakayoko} I. Bakayoko, Hom-post-Lie modules, $\mathcal{O}$-operators and some functors on Hom-algebras, {\it arXiv preprint}, 2016, 	arXiv:1610.02845.


\bibitem{Bai} R. Bai, L. Guo, J. Li, Y. Wu,   Rota-Baxter $3$-Lie algebras, {\it J. Math. Phys.},  {\bf 54},  (2013),  063504.



\bibitem{B60} Baxter, G. An analytic problem whose solution follows from a simple algebraic identity. {\em Pac. J. Math.} {\bf 1960}, {\bf 10}, 731--742.



 \bibitem{Cartan} E. Cartan, Oeuvres completes. Part 1. Gauthier-Villars. Paris. (1952) vol. 2. nos. 101--138.

 \bibitem{Cai}  L. Cai,  Y. Sheng,  Purely Hom-Lie bialgebras,  {\it  Sci. China Math}, {\bf 61}(9), (2018), 1553--1566.



\bibitem{Chen} S. Chen, Q. Lou,  Q. Sun, Cohomologies of Rota-Baxter Lie triple systems and applications, {\it Commun. Algebra}, 2023, https://doi.org/10.1080/00927872.2023.2205938

\bibitem{Chtioui}  T. Chtioui,  A. Hajjaji,  S. Mabrouk,  A. Makhlouf, Cohomologies and deformations of $\mathcal{O}$-operators on Lie triple systems. {\it J. Math. Phys.}, {\bf 64},(2023), 081701.

\bibitem{Das2} A. Das, Cohomology and deformations of weighted Rota-Baxter operators, {\it J. Math. Phys.}, {\bf 63}, (2022),  091703.


\bibitem{Das22}
A. Das, S. Sen,  {Nijenhuis operators on Hom-Lie algebras},
\newblock{\it Commun. Algebra},
\newblock \textbf{50} (2022), 1038--1054.

\bibitem{Das23} A. Das, A. Makhlouf,  {Embedding tensors on Hom-Lie algebras},
\newblock {\it arXiv preprint}, 2023,  arXiv:2304.04178.

\bibitem{Guo0}  L. Guo, W. Keigher, On differential Rota-Baxter algebras, {\it J. Pure Appl. Algebra}, {\bf 212}, (2008), 522--540.




\bibitem{Guo4} S. Guo, Central extensions and deformations of Lie triple systems with a derivation, {\it J.Math. R.Appl.},  {\bf 42}, (2022), 189--198.

\bibitem{Hartwig}  J.Hartwig, D. Larsson, S. Silvestrov, Deformations of Lie algebras using $\sigma$-derivations, {\it  J. Algebra}, {\bf 295}, (2006), 321--344.

\bibitem{Harris} B. Harris,  Cohomology of Lie triple systems and Lie algebras with involution, {\it Trans. Amer. Math. Soc.}, {\bf  98},  (1961), 148--162.

\bibitem{Hou}  S. Hou,  Y. Sheng,  Y. Zhou,  3-post-Lie algebras and relative Rota-Baxter operators of nonzero weight on 3-Lie algebras, {\it J. Algebra},  {\bf 615} (2023), 103--129.

\bibitem{Hou20} Y. Hou, Y. Ma, L. Chen, Product and complex structures on Hom-Lie triple systems. https://www.researchgate.net/publication/341360022.

\bibitem{Hodge} T. Hodge, B.   Parshall, On the representation theory of Lie triple systems, {\it Trans. Amer. Math. Soc.}, {\bf 354}, (2002), 4359--4391.

\bibitem{Hu}  N. Hu, $q$-Witt algebras, $q$-Lie algebras, $q$-holomorph structure and representations, {\it  Algebra Colloq.},  {\bf 6}(1), (1999), 51--70.

\bibitem{Jacobson1} N. Jacobson, Lie and Jordan triple Systems, {\it Amer. J. Math. Soc.}, {\bf 71}(1), { 1949}, 49--170.

\bibitem{Jacobson2} N. Jacobson, General representation theory of Jordan algebras, {\it Trans. Amer. Math. Soc.}, {\bf 70}, (1951), 509--530.

\bibitem{Kubo} F. Kubo,  Y. Taniguchi, A controlling cohomology of the deformation theory of Lie triple systems, {\it J. Algebra},  {\bf 278}, (2004),  242--250.

 \bibitem{K99} Kupershmidt, B.A. What a classical r-matrix really is. {\it J. Nonlinear Math. Phys.} {\bf 1999}, {\em 6}, 448--488.


\bibitem{Li}   Y. Li,  D. Wang,   Relative Rota-Baxter operators on Hom-Lie triple systems, Communications in Algebra, DOI: 10.1080/00927872.2023.2258223

\bibitem{Liu}  S. Liu, L. Song, R. Tang, Representations and cohomologies of regular Hom-pre-Lie algebras, {\it J. Algebra Appl.}, {\bf19},(2020), 2050149 (22 pages).

\bibitem{Lister}  W. Lister,  A structure theory of Lie triple systems, {\it Trans. Amer. Math. Soc.}, {\bf 72}, (1952), 217--242.

\bibitem{Lin} J. Lin,  Y. Wang, S. Deng,  $T^*$-extension of Lie triple systems, {\it Linear Algebra Appl.}, {\bf 431}, (2009), 2071--2083.



\bibitem{Ma} Y. Ma, L. Chen, J. Lin.  Central extensions and deformations of Hom-Lie triple systems, Communications in Algebra, {\bf46}(3),(2018), 1212--1230


 \bibitem{Mishraa18}   S.  Mishra,   A. Naolekar, {$\mathcal{O}$-operators on hom-Lie algebras}, {\it J. Math. Phys.}, \textbf{61} (2020), 121701.




\bibitem{Sun} Q. Sun, S. Chen, Cohomologies and deformations of Lie triple systems with derivations, {\it J. Alg. Appl.}, {\bf 2024}, 2024, 2450053.

\bibitem{Sun17} Q. Sun, H. Li, On parakahler Hom-Lie algebras and hom-left-symmetric bialgebras, Comm. Algebra {\bf45}(1), (2017), 105--120.


\bibitem{Sheng14}  Y. Sheng, C. Bai, {A new approach to Hom-Lie bialgebras}, {\it J. Algebra}, \textbf{399} (2014), 232--250.

\bibitem{Sheng12}  Y. Sheng, {Representations of Hom-Lie algebras}, {\it Alg. Repres. Theo.}, \textbf{15} (2012), 1081--1098.




\bibitem{Teng} W. Teng, F. Long, Y. Zhang, Cohomologies of modified $\lambda$-differential Lie triple systems and applications, {\it AIMS Mathematics}, {\bf8}(10), (2023), 25079--25096.

\bibitem{Teng1} W. Teng, J. Jin, F. Long,  Generalized reynolds operators on Lie-Yamaguti Algebras,  {\it Axioms}, {\bf12}(10), (2023), 934.

\bibitem{Teng2} W. Teng, J. Jin, Y. Zhang, Cohomology of nonabelian embedding tensors on Hom-Lie algebras,   {\it  AIMS Mathematics}, {\bf 8}(9),(2023), 21176--21190.

\bibitem{Teng3} W. Teng, J, Jin,  Relative Rota-Baxter operators on Hom-Lie-Yamaguti algebras, DOI: 10.13140/RG.2.2.33626.00961. https://www.researchgate.net/publication/364144542 (04 October 2022),  appear to  J.Math. R.Appl.



\bibitem{Wang} K. Wang,  G. Zhou, Deformations and homotopy theory of Rota-Baxter algebras of any weight, {\it arXiv preprint}, 2021,  arXiv:2108.06744.



\bibitem{Wu2} X. Wu, Y. Ma, B. Sun,  L. Chen,   Abelian extensions of Lie triple systems with derivations,  {\it Electronic Rese. Arch.}, {\bf 30}(3), (2022), 1087--1103.

\bibitem{Wu22} X. Wu, Y. Ma,  L. Chen,  Relative Rota-Baxter operators of nonzero weights on Lie Triple Systems, {\it arXiv preprint}, 2022,  	arXiv:2207.08946.

\bibitem{Yamaguti} K. Yamaguti, On the cohomology space of Lie triple system, {\it Kumamoto J. Sci. Ser. A}, {\bf 5}, (1960), 44--52.

\bibitem{Yau} D. Yau, On $n$-ary Hom-Nambu and Hom-Nambu-Lie algebras, {\it J. Geom. Phys.}, {\bf 62}(2), (2012),506--522.

\bibitem{Zhou} J. Zhou, L, Chen,  Y.  Ma,  Generalized Derivations of Lie triple systems, {\it Bull.Malaysian Math. Sciences Society}, {\bf 41},  (2018), 637--656.

\bibitem{Zhang} T. Zhang, Notes on cohomologies of Lie triple systems, {\it J. Lie Theorey},  {\bf 24},  (2014),  909--929.

\end{thebibliography}
\end{document}